\documentclass[reqno,12pt]{amsart}
\usepackage{amsmath,amsthm,amssymb,amsfonts,amscd,graphicx}
\input xy
\xyoption{all}
\theoremstyle{plain} 
	\newtheorem{thm}{Theorem}[section]
	\newtheorem*{thm*}{Theorem}
	\newtheorem{cor}[thm]{Corollary}
	\newtheorem{lem}[thm]{Lemma}
	\newtheorem{sublem}[thm]{Sub-Lemma}
	\newtheorem{prop}[thm]{Proposition}
	
	\newtheorem*{conj*}{Conjecture}
	
\theoremstyle{definition}
	\newtheorem{defn}[thm]{Definition}

\theoremstyle{remark}
	\newtheorem{rem}[thm]{Remark}
	
	\newtheorem*{pf}{Proof}
\numberwithin{equation}{section}
\def\CC{{\mathbb C}}

\def\PP{{\mathbb P}}
\def\QQ{{\mathbb Q}}

\def\ZZ{{\mathbb Z}}

\def\F{{\mathcal F}}
\def\G{{\mathcal G}}

\def\I{{\mathcal I}}

\def\M{{\mathcal M}}

\def\O{{\mathcal O}}

\def\T{{\mathcal T}}
\def\X{{\mathcal X}}

\def\p{\partial }

\def\ns{{\nabla}\hspace{-1.4mm}\raisebox{0.3mm}{\text{\footnotesize{\bf /}}}}

\begin{document}
\title[A Uniqueness Theorem for Frobenius Manifolds ]
{A Uniqueness Theorem for Frobenius Manifolds and Gromov--Witten Theory for Orbifold Projective Lines}
\date{\today}
\author{Yoshihisa Ishibashi}
\address{Department of Mathematics, Graduate School of Science, Osaka University, 
Toyonaka Osaka, 560-0043, Japan}
\email{y-ishibashi@cr.math.sci.osaka-u.ac.jp}
\author{Yuuki Shiraishi}
\address{Department of Mathematics, Graduate School of Science, Osaka University, 
Toyonaka Osaka, 560-0043, Japan}
\email{sm5021sy@ecs.cmc.osaka-u.ac.jp}
\author{Atsushi Takahashi}
\address{Department of Mathematics, Graduate School of Science, Osaka University, 
Toyonaka Osaka, 560-0043, Japan}
\email{takahashi@math.sci.osaka-u.ac.jp}
\begin{abstract}
We prove that the Frobenius structure constructed from the Gromov-Witten theory  
for an orbifold projective line with at most three orbifold points is uniquely determined by the WDVV equations with certain 
natural initial conditions. 
\end{abstract}
\maketitle
\section{Introduction}

The Gromov-Witten theory for manifolds is generalized to the one for orbifolds 
by Abramovich--Graber--Vistoli \cite{agv:1} and Chen--Ruan \cite{cr:1}. 
Roughly speaking, their construction of the Gromov-Witten theory for orbifolds is as follows;  
Let $\X$ be an orbifold $($or a smooth proper Deligne--Mumford stack over $\CC)$.
Then, for $g\in\ZZ_{\ge 0}$, $n\in\ZZ_{\ge 0}$ 
and $\beta\in H_2(\X,\ZZ)$, 
the moduli space (stack) $\overline{\M}_{g,n}(\X,\beta)$ 
of orbifold (twisted) stable maps of genus $g$ 
with $n$-marked points of degree $\beta$ is defined.
There exists a virtual fundamental class $[\overline{\M}_{g,n}(\X,\beta)]^{vir}$ 
and Gromov--Witten invariants of genus $g$ 
with $n$-marked points of degree $\beta$ are defined as usual except for 
that we have to use the orbifold cohomology group $H^*_{orb}(\X,\QQ)$:
\[
\left<\Delta_1,\dots, \Delta_n\right>_{g,n,\beta}^\X:=
\int_{[\overline{\M}_{g,n}(\X,\beta)]^{vir}}ev_1^*\Delta_1\wedge \dots \wedge 
ev_n^*\Delta_n,\quad \Delta_1,\dots,\Delta_n\in H^*_{orb}(\X,\QQ),
\]
where $ev^*_i:H^*_{orb}(\X,\QQ)\longrightarrow H^*(\overline{\M}_{g,n}(\X,\beta),\QQ)$
denotes the induced homomorphism by the evaluation map.
We also consider the generating function 
\[
\F_g^\X:=\sum_{n,\beta}\frac{1}{n!}\left<{\bf t},\dots, {\bf t}\right>_{g,n,\beta}^\X,
\quad {\bf t}=\sum_{i}t_i\Delta_i
\]
and call it the genus $g$ potential where $\{\Delta_i\}$ denotes a $\QQ$-basis of 
$H^*_{orb}(\X,\QQ)$.
The main result in \cite{agv:1, cr:1} tells us that 
we can treat the Gromov--Witten theory defined for orbifolds 
as if $\X$ were a usual manifold.
In particular, we have the point axiom,
the divisor axiom for a class in $H^2(\X,\QQ)$
and the associativity of the quantum product, namely, the Witten--Dijkgraaf--Verlinde--Verlinde $($WDVV$)$ equations
(see \cite{agv:1, cr:1} for details of these axioms.), 
which gives a formal Frobenius manifold $M$ whose structure sheaf $\O_M$, tangent sheaf $\T_M$ and
Frobenius potential are defined as 
the algebra $\CC[[H^*_{orb}(\X,\CC)]]$ of formal power series in dual coordinates $\{t_i\}$ 
of the $\QQ$-basis $\{\Delta_i\}$ of $H^*_{orb}(\X,\QQ)$,  
$\T_M:=H^*_{orb}(\X,\CC)\otimes_\CC\CC[[H^*_{orb}(\X,\CC)]]$ and the genus zero 
potential $\F_0^\X$, respectively.
Let $A$ be a triplet $(a_1,a_2,a_3)$ of positive integers such that $a_1\le a_2\le a_3$.
Set $\mu_A=a_1+a_2+a_3-1$ and $\chi_A:=1/a_1+1/a_2+1/a_3-1$.
We shall consider the orbifold projective line with three orbifold points whose orders are 
$a_1, a_2, a_3$, which is denoted by $\PP^{1}_A$.
There, $\mu_A$ is regarded as the total dimension of the orbifold cohomology group $H^*_{orb}(\X,\CC)$ and  
$\chi_A$ is regarded as the orbifold Euler number of $\PP^1_A$. 
In this paper we show that the Frobenius structure constructed from the Gromov--Witten theory for $\PP^{1}_A$ 
can be determined by the WDVV equations with certain natural initial conditions.  
First, we shall show the following uniqueness theorem, which is our main theorem in this paper:
\begin{thm*}[Theorem~\ref{first}]
There exists a unique Frobenius manifold $M$ of rank $\mu_A$ and dimension one with flat coordinates 
$(t_1,t_{1,1},\dots ,t_{i,j},\dots ,t_{3,a_3-1},t_{\mu_A})$ satisfying the following conditions$:$
\begin{enumerate}
\item 
The unit vector field $e$ and the Euler vector field $E$ are given by
\[
e=\frac{\p}{\p t_1},\ E=t_1\frac{\p}{\p t_1}+\sum_{i=1}^3\sum_{j=1}^{a_i-1}\frac{a_i-j}{a_i}t_{i,j}\frac{\p}{\p t_{i,j}}
+\chi_A\frac{\p}{\p t_{\mu_A}}.
\]
\item 
The non-degenerate symmetric bilinear form $\eta$ on $\T_M$ satisfies
\begin{align*}
&\ \eta\left(\frac{\p}{\p t_1}, \frac{\p}{\p t_{\mu_A}}\right)=
\eta\left(\frac{\p}{\p t_{\mu_A}}, \frac{\p}{\p t_1}\right)=1,\\ 
&\ \eta\left(\frac{\p}{\p t_{i_1,j_1}}, \frac{\p}{\p t_{i_2,j_2}}\right)=
\begin{cases}
\frac{1}{a_{i_1}}\quad i_1=i_2\text{ and }j_2=a_{i_1}-j_1,\\
0 \quad \text{otherwise}.
\end{cases}
\end{align*}
\item 
The Frobenius potential $\F$ satisfies $E\F|_{t_{1}=0}=2\F|_{t_{1}=0}$,
\[
\left.\F\right|_{t_1=0}\in\CC\left[[t_{1,1}, \dots, t_{1,a_1-1},
t_{2,1}, \dots, t_{2,a_2-1},t_{3,1}, \dots, t_{3,a_3-1},e^{t_{\mu_A}}]\right].
\]
\item Assume the condition {\rm (iii)}. we have
\begin{equation*}
\F|_{t_1=e^{t_{\mu_A}}=0}=\G^{(1)}+\G^{(2)}+\G^{(3)}, \quad \G^{(i)}\in \CC[[t_{i,1},\dots, t_{i,a_i-1}]],\ i=1,2,3.
\end{equation*}
\item 
Assume the condition {\rm (iii)}. In the frame $\frac{\p}{\p t_1}, \frac{\p}{\p t_{1,1}},\dots, 
\frac{\p}{\p t_{3,a_3-1}},\frac{\p}{\p t_{\mu_A}}$ of $\T_M$,
the product $\circ$ can be extended to the limit $t_1=t_{1,1}=\dots=t_{3,a_3-1}=e^{t_{\mu_A}}=0$.
The $\CC$-algebra obtained in this limit is isomorphic to
\[
\CC[x_1,x_2,x_3]\left/\left(x_1x_2,x_2x_3,x_3x_1,a_1x_1^{a_1}-a_2x_2^{a_2},
a_2x_2^{a_2}-a_3x_3^{a_3},a_3x_3^{a_3}-a_1x_1^{a_1}\right)\right.,
\]
where $\p/\p t_{1,1},\p/\p t_{2,1},\p/\p t_{3,1}$ are mapped to
$x_1,x_2,x_3$, respectively.
\item The term 
\[
\begin{cases}
e^{t_{\mu_A}}\quad \textit{if}\quad a_{1}=a_{2}=a_{3}=1,\\
t_{3,1}e^{t_{\mu_A}}\quad \textit{if}\quad 1=a_{1}=a_{2}<a_{3},\\
t_{2,1}t_{3,1}e^{t_{\mu_A}}\quad \textit{if}\quad 1=a_{1}<a_{2},\\
t_{1,1}t_{2,1}t_{3,1}e^{t_{\mu_A}}\quad \textit{if}\quad a_{1}\ge 2,
\end{cases}
\]
occurs with the coefficient $1$ in $\F$. 
\end{enumerate}
\end{thm*}
Actually, the condition {\rm (iv)} follows from others if $A\ne (1,2,2)$ or $A\ne (2,2,a_3)$, $a_3\ge 2$. 
If $A=(1,2,2)$ or $A=(2,2,a_3)$, $a_3\ge 3$, then we only have to assume instead of the condition {\rm (iv)} 
the weaker and more natural one 
\begin{itemize}
\item[{\rm (iv')}]
If $a_{i_1}=a_{i_2}$ for some $i_1,i_2\in\{1,2,3\}$, then the
Frobenius potential $\F$ is invariant under
the permutation of parameters $t_{i_1,j}$ and $t_{i_2,j}$ $(j=1,\dots, a_{i_1}-1)$.
\end{itemize}
Therefore, in the proof of Theorem~\ref{first}, 
we shall also show the condition {\rm (iv)} in Proposition~\ref{sep} (under the condition {\rm (iv')} when $A=(2,2,a_3)$, $a_3\ge 3$
or $A=(1,2,2)$).
Next, we show the following.
\begin{thm*}[Theorem~\ref{second}]
The conditions in Theorem~\ref{first} are satisfied by the Frobenius structure constructed
from the Gromov--Witten theory for $\PP^1_A$. 
\end{thm*}
As a corollary, the Frobenius structure constructed from the Gromov--Witten theory for $\PP^{1}_A$ can be uniquely 
reconstructed by the conditions in Theorem~\ref{first}.  
In our subsequent paper, we shall show that the Frobenius structure constructed from a universal unfolding of the polynomial 
$f_A(x_1,x_2,x_3):=x_1^{a_1}+x_2^{a_2}+x_3^{a_3}-q^{-1}x_1x_2x_3$, $q\in\CC\backslash\{0\}$ 
via a fixed primitive form $\zeta_A$ also satisfies the conditions in Theorem~\ref{first}.
As a consequence, we obtain the mirror isomorphism as Frobenius manifolds 
between the one constructed from the Gromov--Witten theory 
for $\PP^{1}_A$ and the one constructed from the pair $(f_A,\zeta_A)$. This result is already obtained by  
Milanov--Tseng \cite{mt:1} for the case $a_1=1$ and by Rossi \cite{r:1} for the case $\chi_A>0$. The present work
will be the first in a series of papers to simplify and generalise their works. 
\bigskip
\noindent
{\it Acknowledgement}\\
\indent
The second named author is deeply grateful to Claus Hertling and 
Christian Sevenheck for their valuable discussions and encouragement.  
He is supported by the JSPS International Training Program (ITP).
The third named author is supported by Grant-in Aid for Scientific Research 
grant numbers 24684005 from the Ministry of Education, 
Culture, Sports, Science and Technology, Japan.


\section{Preliminary}

In this section, we recall the definition and some basic properties of the Frobenius manifold \cite{d:1}. 
The definition below is taken from Saito-Takahashi \cite{st:1}.
\begin{defn}
Let $M=(M,\O_{M})$ be a connected complex manifold or a formal manifold over $\CC$ of dimension $\mu$
whose holomorphic tangent sheaf and cotangent sheaf are denoted by $\T_{M}, \Omega_M^1$ respectively
and $d$ a complex number.
A {\it Frobenius structure of
rank $\mu$ and dimension $d$ on M} is a tuple $(\eta, \circ , e,E)$, where $\eta$ is a non-degenerate $\O_{M}$-symmetric bilinear
form on $\T_{M}$, $\circ $ is $\O_{M}$-bilinear product on $\T_{M}$, defining an associative and commutative
$\O_{M}$-algebra structure with the unit $e$, and $E$ is a holomorphic vector field on $M$, called
the Euler vector field, which are subject to the following axioms:
\begin{enumerate}
\item The product $\circ$ is self-ajoint with respect to $\eta$: that is,
\begin{equation*}
\eta(\delta\circ\delta',\delta'')=\eta(\delta,\delta'\circ\delta''),\quad
\delta,\delta',\delta''\in\T_M. 
\end{equation*} 
\item The {\rm Levi}--{\rm Civita} connection $\ns:\T_M\otimes_{\O_M}\T_M\to\T_M$ with respect to $\eta$ is
flat: that is, 
\begin{equation*}
[\ns_\delta,\ns_{\delta'}]=\ns_{[\delta,\delta']},\quad \delta,\delta'\in\T_M.
\end{equation*}
\item The tensor $C:\T_M\otimes_{\O_M}\T_M\to \T_M$  defined by 
$C_\delta\delta':=\delta\circ\delta'$, $(\delta,\delta'\in\T_M)$ is flat: that is,
 \begin{equation*}
\ns C=0.
\end{equation*} 
\item The unit element $e$ of the $\circ $-algebra is a $\ns$-flat homolophic vector field: that is,
\begin{equation*}
\ns e=0.
\end{equation*} 
\item The metric $\eta$ and the product $\circ$ are homogeneous of degree $2-d$ ($d\in\CC$), $1$ respectively 
with respect to Lie derivative $Lie_{E}$ of {\rm Euler} vector field $E$: that is,
\begin{equation*}
Lie_E(\eta)=(2-d)\eta,\quad Lie_E(\circ)=\circ.
\end{equation*}
\end{enumerate}
\end{defn}
Next we expose some basic properties of the Frobenius manifold without their proofs.
Let us consider the space of horizontal sections of the connection $\ns$:
\[
\T_M^f:=\{\delta\in\T_M~|~\ns_{\delta'}\delta=0\text{ for all }\delta'\in\T_M\}
\]
which is a local system of rank $\mu $ on $M$ such that the metric $\eta$ 
takes constant value on $\T_M^f$. Namely, we have 
\begin{equation}
\eta (\delta,\delta')\in\CC,\quad  \delta,\delta' \in \T_M^f.
\end{equation}
\begin{prop}\label{prop:flat coordinates}
At each point of $M$, there exist a local coordinate $(t_1,\dots,t_{\mu})$, called flat coordinates, such that
$e=\p_1$, $\T_M^f$ is spanned by $\p_1,\dots, \p_{\mu}$ and $\eta(\p_i,\p_j)\in\CC$ for all $i,j=1,\dots, \mu$,
where we denote $\p/\p t_i$ by $\p_i$. 
\end{prop} 
The axiom $\ns C=0$, implies the following:
\begin{prop}\label{prop:potential}
At each point of $M$, there exist the local holomorphic function $\F$, called Frobenius potential, satisfying
\begin{equation*}
\eta(\p_i\circ\p_j,\p_k)=\eta(\p_i,\p_j\circ\p_k)=\p_i\p_j\p_k \F,
\quad i,j,k=1,\dots,\mu,
\end{equation*}
for any system of flat coordinates. In particular, one has
\begin{equation*}
\eta_{ij}:=\eta(\p_i,\p_j)=\p_1\p_i\p_j \F. 
\end{equation*}
\end{prop}
\section{A Uniqueness Theorem}

We have the following uniqueness theorem for Frobenius manifolds of rank $\mu_A$ and dimension one.

\begin{thm}\label{first}
There exists a unique Frobenius manifold $M$ of rank $\mu_A$ and dimension one with flat coordinates 
$(t_1,t_{1,1},\dots ,t_{i,j},\dots ,t_{3,a_3-1},t_{\mu_A})$ satisfying the following conditions$:$
\begin{enumerate}
\item 
The unit vector field $e$ and the Euler vector field $E$ are given by
\[
e=\frac{\p}{\p t_1},\ E=t_1\frac{\p}{\p t_1}+\sum_{i=1}^3\sum_{j=1}^{a_i-1}\frac{a_i-j}{a_i}t_{i,j}\frac{\p}{\p t_{i,j}}
+\chi_A\frac{\p}{\p t_{\mu_A}}.
\]
\item 
The non-degenerate symmetric bilinear form $\eta$ on $\T_M$ satisfies
\begin{align*}
&\ \eta\left(\frac{\p}{\p t_1}, \frac{\p}{\p t_{\mu_A}}\right)=
\eta\left(\frac{\p}{\p t_{\mu_A}}, \frac{\p}{\p t_1}\right)=1,\\ 
&\ \eta\left(\frac{\p}{\p t_{i_1,j_1}}, \frac{\p}{\p t_{i_2,j_2}}\right)=
\begin{cases}
\frac{1}{a_{i_1}}\quad i_1=i_2\text{ and }j_2=a_{i_1}-j_1,\\
0 \quad \text{otherwise}.
\end{cases}
\end{align*}
\item 
The Frobenius potential $\F$ satisfies $E\F|_{t_{1}=0}=2\F|_{t_{1}=0}$,
\[
\left.\F\right|_{t_1=0}\in\CC\left[[t_{1,1}, \dots, t_{1,a_1-1},
t_{2,1}, \dots, t_{2,a_2-1},t_{3,1}, \dots, t_{3,a_3-1},e^{t_{\mu_A}}]\right]
\]
\item Assume the condition {\rm (iii)}. we have
\begin{equation*}
\F|_{t_1=e^{t_{\mu_A}}=0}=\G^{(1)}+\G^{(2)}+\G^{(3)}, \quad \G^{(i)}\in \CC[[t_{i,1},\dots, t_{i,a_i-1}]],\ i=1,2,3.
\end{equation*}
\item 
Assume the condition {\rm (iii)}. In the frame $\frac{\p}{\p t_1}, \frac{\p}{\p t_{1,1}},\dots, 
\frac{\p}{\p t_{3,a_3-1}},\frac{\p}{\p t_{\mu_A}}$ of $\T_M$,
the product $\circ$ can be extended to the limit $t_1=t_{1,1}=\dots=t_{3,a_3-1}=e^{t_{\mu_A}}=0$.
The $\CC$-algebra obtained in this limit is isomorphic to
\[
\CC[x_1,x_2,x_3]\left/\left(x_1x_2,x_2x_3,x_3x_1,a_1x_1^{a_1}-a_2x_2^{a_2},
a_2x_2^{a_2}-a_3x_3^{a_3},a_3x_3^{a_3}-a_1x_1^{a_1}\right)\right.,
\]
where $\p/\p t_{1,1},\p/\p t_{2,1},\p/\p t_{3,1}$ are mapped to
$x_1,x_2,x_3$, respectively.
\item The term 
\[
\begin{cases}
e^{t_{\mu_A}}\quad \textit{if}\quad a_{1}=a_{2}=a_{3}=1,\\
t_{3,1}e^{t_{\mu_A}}\quad \textit{if}\quad 1=a_{1}=a_{2}<a_{3},\\
t_{2,1}t_{3,1}e^{t_{\mu_A}}\quad \textit{if}\quad 1=a_{1}<a_{2},\\
t_{1,1}t_{2,1}t_{3,1}e^{t_{\mu_A}}\quad \textit{if}\quad a_{1}\ge 2,
\end{cases}
\]
occurs with the coefficient $1$ in $\F$. 
\end{enumerate}
\end{thm}
\begin{rem}
As is already explained in Introduction, in the proof of Theorem~\ref{first}, 
we shall also show in Proposition~\ref{sep} that the condition {\rm (iv)} follows from others 
if $a_3\ge 3$ and if $A=(1,2,2)$ or $A=(2,2,a_3)$, $a_3\ge 3$ together with the following condition {\rm (iv')}$;$
\begin{itemize}
\item[{\rm (iv')}]
If $a_{i_1}=a_{i_2}$ for some $i_1,i_2\in\{1,2,3\}$, then the
Frobenius potential $\F$ is invariant under
the permutation of parameters $t_{i_1,j}$ and $t_{i_2,j}$ $(j=1,\dots, a_{i_1}-1)$.
\end{itemize}
Namely, only for $A=(2,2,2)$, the stronger condition {\rm (iv)} is necessary.
\end{rem}

By the condition {\rm (iii)} in Theorem~\ref{first},
we can expand the non-trivial part of the Frobenius potential $\F|_{t_1=0}$ as
\[
\F|_{t_{1}=0}=\sum_{\alpha =(\alpha_{1,1},\dots ,\alpha_{3,a_{3}-1})} c(\alpha ,m) t^{\alpha}e^{mt_{\mu_A}}, \ \
t^{\alpha}=\prod_{i=1}^3\prod_{j=1}^{a_i-1}t_{i,j}^{\alpha_{i,j}}.
\]
Consider a free abelian group $\ZZ^{\mu_A-2}$ and denote  its 
standard basis by $e_{i,j}$, $i=1,2,3$, $j=1,\dots, a_i-1$.
The element $\alpha=\sum_{i,j}\alpha_{i,j}e_{i,j}$, $\alpha_{i,j}\in\ZZ$ of  $\ZZ^{\mu_A-2}$ is called {\it non-negative} 
and is denoted by  $\alpha\ge 0$ if all $\alpha_{i,j}$ are non-negative integers.
We also denote by $c(e_{1}+e_{i,j}+e_{i,a_{i}-j},0)$ the coefficient of $t_{1}t_{i,j}t_{i,a_{i}-j}$ in the trivial part of the Frobenius potential $\F$.
For a non-negative $\alpha\in\ZZ^{\mu_A-2}$, we set 
\[
|\alpha|:=\sum_{i=1}^3\sum_{j=1}^{a_i-1}\alpha_{i,j},
\]
and call it the {\it length} of $\alpha$. 
Define the number $s_{a,b,c}$ for $a,b,c\in \ZZ$ as follows:
\[
s_{a,b,c}=
\begin{cases}
1 \ \ \text{if} \ \ a,b,c \text{ are distinct},\\
6 \ \ \text{if} \ \ a=b=c,\\
2 \ \ \text{otherwise}.\\
\end{cases}
\]
For $a, b, c, d\in \{1,\dots, \mu_A\}$, denote by $WDVV(a,b,c,d)$   
the following equation:
\[
\displaystyle \sum_{\sigma ,\tau=1}^{\mu_A}\p_{a}\p_{b}\p_{\sigma}\F \cdot \eta^{\sigma \tau}\cdot \p_{\tau}\p_{c}\p_{d}\F
-\sum_{\sigma ,\tau=1}^{\mu_A}\p_{a}\p_{c}\p_{\sigma}\F \cdot \eta^{\sigma \tau}\cdot \p_{\tau}\p_{b}\p_{d}\F=0,
\]
where $(\eta^{\sigma \tau}):=(\eta_{\sigma \tau})^{-1}$.


\subsection{$c(\alpha,0)$ and $c(\alpha,1)$ can be reconstructed}\label{reconst-s1}
\begin{prop}\label{lem3}
$c(\alpha,0)$ with $|\alpha|=3$ are determined by the condition {\rm (v)} of Theorem~\ref{first}.
\end{prop}

\begin{pf}
Note that $C_{ijk}=\eta (\partial_{i}\circ \partial_{j}, \partial_{k})$
and the metric $\eta$ can be extended to the limit $\underline{t},e^{t} \rightarrow 0$.
We denote this extended metric by $\eta'$. 
By the condition {\rm (v)}, the relation $xy=yz=zx=0$ holds in the $\CC$-algebra obtained in this limit.
Therefore, $c(\sum_{k=1}^3e_{i_k,j_{k}},0)\ne 0$ only if $i_1=i_2=i_3$.
One has 
\begin{multline*}
s_{j_{1},j_{2},j_{3}}\cdot c\left(\sum_{k=1}^3e_{i,j_{k}},0\right)=\lim_{\underline{t},e^{t} \rightarrow 0} \p_{i,j_{1}}\p_{i,j_{2}}\p_{i,j_{3}}\F
=\eta' (x_{i}^{j_{1}}\cdot x_{i}^{j_{2}},x_{i}^{j_{3}})\\
=\eta' (1\cdot x_{i}^{j_1+j_2}, x_{i}^{j_{3}})
=\lim_{\underline{t},e^{t} \rightarrow 0} \p_{1}\p_{i, j_{1}+j_{2}}\p_{i,j_{3}}\F
\end{multline*}
and 
\[
\lim_{\underline{t},e^{t} \rightarrow 0} \p_{1}\p_{i, j_{1}+j_{2}}\p_{i, j_{3}}\F=
\begin{cases}
\frac{1}{a_i}\quad \text{if } \sum_{k=1}^3j_{k}=a_i,\\
0\quad \text{otherwise}.
\end{cases}
\]
\qed
\end{pf}

The following Lemma~\ref{wsep} is a part of the condition {\rm (iv)} in Theorem~\ref{first}.
However, for the later convenience, we shall show that Lemma~\ref{wsep} is derived from other conditions in Theorem~\ref{first}.
\begin{lem}\label{wsep}
Let $\gamma \in \ZZ^{\mu_{A}-2}$ be a non-negative element such that 
$|\gamma|=4$ and $\gamma -e_{i_{1},j_{1}}-e_{i_{2},j_{2}}\ge 0$ for $i_{1}\ne i_{2}$.
If $a_{i_1}\ge 3$, then one has $c(\gamma,0)=0$.
\end{lem}
\begin{pf}
Note that $c(\beta,0)=0$ if $|\beta|=3$ and $\beta-e_{i_{1},j_{1}}-e_{i_{2},j_{2}}\ge 0$ for $i_{1}\ne i_{2}$ by Proposition \ref{lem3}.
We shall split the proof into the following two steps:

\vspace{5pt}
\noindent
Step 1: We shall consider the case that the term $t^{\gamma}$ has, as a factor,
$t_{i_1,j}$ such that $a_{i_1}\ge 3$ and $j\ge 2$ for some $j$.
We shall split, moreover, Step 1 into following four cases: 
\vspace{5pt}
\begin{enumerate}
\item The term $t^{\gamma}$ has $t_{1,j_1}, t_{2,j_2}, t_{3,j_3}$ as factors.
\item The term $t^{\gamma}$ has $t_{i,j}$ and $t_{i,j'}$ as factors for each $i=i_{1}, i_{2}$.
\item The term $t^{\gamma}$ has $t_{i_1,j_1}, t_{i_1,j'_1}, t_{i_1,j''_1}$ and only $t_{i_2,j_2}$ as factors.
\item The term $t^{\gamma}$ has, as factors, $t_{i_1,j_1}$ and $t_{i_{2},j}, t_{i_{2},j'}, t_{i_{2},j''}$.
\end{enumerate}
\vspace{5pt}
\noindent
Step 2: We shall consider the case that $t^{\gamma}$ does not have, as factors, $t_{i,j}$ for $j\ge 2$. 

\vspace{5pt}
\noindent
\underline{\it Step 1-{\rm (i)}. The term $t^{\gamma}$ has $t_{1,j_1}, t_{2,j_2}, t_{3,j_3}$ as factors. }

\begin{sublem}[Step 1-(i)]\label{sl1-l1-wsep}
Let $\gamma \in \ZZ^{\mu_{A}-2}$ be a non-negative element such that $|\gamma|=4$ and $\gamma -e_{i_{1},j_{1}}-e_{i_{2},j_{2}}-e_{i_{3},j_{3}}\ge 0$ for distinct $i_{1}$, $i_{2}$, $i_{3}$.
If $a_{i_1}\ge 3$ and $j_{1}\ge 2$, then one has $c(\gamma,0)=0$. 
\end{sublem}

\begin{pf}
Taking the coefficient in front of $t^{\gamma -e_{i_{1},j_{1}}-e_{i_{2},j_{2}}-e_{i_{3},j_{3}}}$ 
in $WDVV((i_{3},j_{3}),(i_{2},j_{2}),(i_{1},j_{1}-1),(i_{1},1))$,
one has
\[
\gamma_{i_{1},j_{1}}\gamma_{i_{2},j_{2}}\gamma_{i_{3},j_{3}}\cdot c(\gamma,0)\cdot a_{i_{1}}\cdot 
s_{a_{i_{1}}-j_{1},j_{1}-1,1}\cdot c(e_{i_{1},a_{i_{1}}-j_{1}}+e_{i_{1},j_{1}-1}+e_{i_{1},1},0)=0.
\]
Hence one has $c(\gamma ,0)=0$.
\qed
\end{pf}

\vspace{5pt}
\noindent
\underline{\it Step 1-{\rm (ii)}.  The term
$t^{\gamma}$ has $t_{i,j}$ and $t_{i,j'}$ as factors for each $i=i_{1}, i_{2}$.}
 
\begin{sublem}[Step 1-(ii)]\label{sl1-l2-wsep}
Let
$\gamma \in \ZZ^{\mu_A-2}$ be $e_{i_{1},j_{1}}+e_{i_{1},j'_{1}}+e_{i_{2},j_{2}}+e_{i_{2},j'_{2}}$ for $i_{1}\neq i_{2}$. 
If $a_{i_{1}}\geq 3$ and $j'_{1}\ge 2$, then one has $c(\gamma,0)=0$.
\end{sublem}

\begin{pf}
Taking the coefficient in front of $t_{i_2,j'_2}$ in $WDVV((i_{2},j_{2}),(i_{1},j_{1}),(i_{1},1),(i_{1},j'_{1}-1))$, 
one has
\begin{eqnarray*}
\lefteqn{{\rm (i)} \ \gamma_{i_{1},j_{1}}\gamma_{i_{2},j_{2}}\gamma_{i_{1},j'_{1}}\cdot c(\gamma,0)\cdot a_{i_{1}}\cdot 
s_{1,j'_{1}-1,a_{i_{1}}-j'_{1}}\cdot c(e_{i_{1},1}+e_{i_{1},j'_{1}-1}+e_{i_{1},a_{i_{1}}-j'_{1}},0)}\\
&&-(\gamma'_{i_{1},1}+1)(\gamma'_{i_{2},j_{2}}+1)(\gamma'_{i_{1},j_{1}+j'_{1}-1}+1)\cdot c(\gamma' +e_{i_{1},1}+e_{i_{2},j_{2}}+e_{i_{1},j_{1}+j'_{1}-1},0)
\cdot a_{i_{1}}\cdot \\
&&s_{j_1,j'_{1}-1,a_{i_{1}}+1-j_{1}-j'_{1}}\cdot c(e_{i_{1},j_{1}}+e_{i_{1},j'_{1}-1}+e_{i_{1},a_{i_{1}}+1-j_{1}-j'_{1}},0)=0\\ 
&&\text{if} \ \ 3\leq j_{1}+j'_{1}\leq a_{i_{1}} \ \text{and} \ \text{where} \ \gamma'=\gamma-e_{i_{1},j_{1}}-e_{i_{2},j_{2}}-e_{i_{1},j'_{1}}=e_{i_2,j'_2},\\ 
\lefteqn{{\rm (ii)} \ \gamma_{i_{1},j_{1}}\gamma_{i_{2},j_{2}}\gamma_{i_{1},j'_{1}}\cdot c(\gamma,0)\cdot a_{i_{1}}\cdot s_{1,j'_{1}-1,a_{i_{1}}-j'_{1}}\cdot
c(e_{i_{1},1}+e_{i_{1},j'_{1}-1}+e_{i_{1},a_{i_1}-j'_{1}},0)=0}\\
&&\text{if} \ \ j_{1}+j'_{1}> a_{i_{1}}.
\end{eqnarray*}
We shall show that $c(\gamma' +e_{i_{1},1}+e_{i_{2},j_{2}}+e_{i_{1},j_{1}+j'_{1}-1},0)=0$. One has $\deg(t_{i_{2},j_{2}}t_{i_{2},j'_{2}})=(j_{1}+j'_{1})/a_{i_{1}}\le 1$. 
Taking the coefficient in front of $t_{i_{1},1}$ in $WDVV((i_{2},j_{2}),(i_{2},j'_{2}),(i_{1},j_{1}+j'_{1}-2),(i_{1},1))$, 
one has
\begin{multline*}
(\gamma'_{i_{2},j_{2}}+1)(\gamma'_{i_{2},j'_{2}}+1)(\gamma'_{i_{1},j_{1}+j'_{1}-1}+1)\cdot c(\gamma' +e_{i_{1},1}+e_{i_{2},j_{2}}+e_{i_{1},j_{1}+j'_{1}-1},0)
\cdot a_{i_{1}} \cdot \\
s_{1,j_{1}+j'_{1}-2,a_{i_{1}}+1-j_{1}-j'_{1}}\cdot c(e_{i_{1},1}+e_{i_{1},j_{1}+j'_{1}-2}+e_{i_{1},a_{i_{1}}+1-j_{1}-j'_{1}},0)=0.
\end{multline*}
By Proposition \ref{lem3}, $c(\gamma' +e_{i_{1},1}+e_{i_{2},j_{2}}+e_{i_{1},j_{1}+j'_{1}-1},0)=0$.
Hence one has $c(\gamma,0)=0$.
\qed
\end{pf}

\vspace{5pt}
\noindent
\underline{\it Step 1-{\rm (iii)}. The term
$t^{\gamma}$ has $t_{i_1,j_1}, t_{i_1,j'_1}, t_{i_1,j''_1}$ and only $t_{i_2,j_2}$ as factors.}

\vspace{10pt}
We shall split Step 1-{\rm (iii)} into Case 1 ($j_2\ge 2$) and Case 2 ($j_2=1$).
\vspace{-5pt}
\begin{sublem}[Step 1-(iii)-Case 1]\label{sl2-l2-wsep}
Let 
$\gamma \in \ZZ^{\mu_{A}-2}$ be 
$e_{i_{1},j_{1}}+e_{i_{1},j'_{1}}+e_{i_{1},j''_{1}}+e_{i_{2},j_{2}}$ for $i_{1}\neq i_{2}$. 
If $a_{i_{1}}\ge 3$, $a_{i_{2}}\ge 3$, $j'_{1}\ge 2$ and $j_{2}\ge 2$, then one has $c(\gamma,0)=0$.
\end{sublem}

\begin{pf}
Taking the coefficient in front of $t_{i_1,j''_1}$ in $WDVV((i_{2},j_{2}),(i_{1},j_{1}),(i_{1},1),(i_{1},j'_{1}-1))$, 
one has
\begin{eqnarray*}
\lefteqn{{\rm (i)} \ \gamma_{i_{1},j_{1}}\gamma_{i_{2},j_{2}}\gamma_{i_{1},j'_{1}}\cdot c(\gamma,0)\cdot a_{i_{1}}\cdot 
s_{1,j'_{1}-1,a_{i_{1}}-j'_{1}}\cdot c(e_{i_{1},1}+e_{i_{1},j'_{1}-1}+e_{i_{1},a_{i_{1}}-j'_{1}},0)}\\
&&-(\gamma'_{i_{1},1}+1)(\gamma'_{i_{2},j_{2}}+1)(\gamma'_{i_{1},j_{1}+j'_{1}-1}+1)\cdot c(\gamma' +e_{i_{1},1}+e_{i_{2},j_{2}}+e_{i_{1},j_{1}+j'_{1}-1},0)
\cdot a_{i_{1}}\cdot \\
&&s_{j_1,j'_{1}-1,a_{i_{1}}+1-j_{1}-j'_{1}}\cdot c(e_{i_{1},j_{1}}+e_{i_{1},j'_{1}-1}+e_{i_{1},a_{i_{1}}+1-j_{1}-j'_{1}},0)=0\\ 
&&\text{if} \ \ 3\leq j_{1}+j'_{1}\leq a_{i_{1}} \ \text{and} \ \text{where} \ \gamma'=\gamma-e_{i_{1},j_{1}}-e_{i_{2},j_{2}}-e_{i_{1},j'_{1}}=e_{i_1,j''_1},\\ 
\lefteqn{{\rm (ii)} \ \gamma_{i_{1},j_{1}}\gamma_{i_{2},j_{2}}\gamma_{i_{1},j'_{1}}\cdot c(\gamma,0)\cdot a_{i_{1}}\cdot s_{1,j'_{1}-1,a_{i_{1}}-j'_{1}}\cdot
c(e_{i_{1},1}+e_{i_{1},j'_{1}-1}+e_{i_{1},a_{i_1}-j'_{1}},0)=0}\\
&&\text{if} \ \ j_{1}+j'_{1}> a_{i_{1}}.
\end{eqnarray*}
We shall show that $c(\gamma'+e_{i_{1},1}+e_{i_{2},j_{2}}+e_{i_{1},j_{1}+j'_{1}-1},0)=0$.
One has $\deg(t_{i_{1},j''_{1}})\leq (j_{1}+j'_{1}-1)/a_{i_{1}}$. Then taking 
the coefficient in front of $t_{i_{1},1}$ in $WDVV((i_{1},j''_{1}),(i_{1},j_{1}+j'_{1}-1),(i_{2},j_{2}-1),(i_{2},1))$,
one has 
\begin{multline*}
\lefteqn{(\gamma'_{i_{1},j''_{1}}+1)(\gamma'_{i_{1},j_{1}+j'_{1}-1}+1)(\gamma'_{i_{2},j_{2}}+1)\cdot c(\gamma' +e_{i_{1},1}+e_{i_{2},j_{2}}+e_{i_{1},j_{1}+j'_{1}-1},0)
\cdot a_{i_{2}}\cdot}\\
s_{1,j_{2}-1,a_{i_{2}}-j_{2}}\cdot c(e_{i_{2},1}+e_{i_{2},j_{2}-1}+e_{i_{2},a_{i_{2}}-j_{2}},0)=0.
\end{multline*}
By Proposition \ref{lem3}, $c(\gamma'+e_{i_{1},1}+e_{i_{2},j_{2}}+e_{i_{1},j_{1}+j'_{1}-1},0)=0$.
Hence one has $c(\gamma,0)=0$.
\qed
\end{pf}

\begin{sublem}[Step 1-(iii)-Case 2]\label{sl3-l2-wsep}
Let
$\gamma \in \ZZ^{\mu_{A}-2}$ be 
$e_{i_{1},j_{1}}+e_{i_{1},j'_{1}}+e_{i_{1},j''_{1}}+e_{i_{2},j_{2}}$ for $i_{1}\neq i_{2}$. 
If $a_{i_{1}}\ge 3$ $j'_{1}\ge 2$ and $j_{2}=1$, then one has $c(\gamma,0)=0$.
\end{sublem}

\begin{pf}
Taking the coefficient in front of $t_{i_1,j''_1}$ in $WDVV((i_{2},j_{2}),(i_{1},j_{1}),(i_{1},1),(i_{1},j'_{1}-1))$, 
one has
\begin{eqnarray*}
\lefteqn{{\rm (i)} \ \gamma_{i_{1},j_{1}}\gamma_{i_{2},j_{2}}\gamma_{i_{1},j'_{1}}\cdot c(\gamma,0)\cdot a_{i_{1}}\cdot 
s_{1,j'_{1}-1,a_{i_{1}}-j'_{1}}\cdot c(e_{i_{1},1}+e_{i_{1},j'_{1}-1}+e_{i_{1},a_{i_{1}}-j'_{1}},0)}\\
&&-(\gamma'_{i_{1},1}+1)(\gamma'_{i_{2},j_{2}}+1)(\gamma'_{i_{1},j_{1}+j'_{1}-1}+1)\cdot c(\gamma' +e_{i_{1},1}+e_{i_{2},j_{2}}+e_{i_{1},j_{1}+j'_{1}-1},0)
\cdot a_{i_{1}}\cdot \\
&&s_{j_1,j'_{1}-1,a_{i_{1}}+1-j_{1}-j'_{1}}\cdot c(e_{i_{1},j_{1}}+e_{i_{1},j'_{1}-1}+e_{i_{1},a_{i_{1}}+1-j_{1}-j'_{1}},0)=0\\ 
&&\text{if} \ \ 3\leq j_{1}+j'_{1}\leq a_{i_{1}} \ \text{and} \ \text{where} \ \gamma'=\gamma-e_{i_{1},j_{1}}-e_{i_{2},j_{2}}-e_{i_{1},j'_{1}}=e_{i_1,j''_1},\\ 
\lefteqn{{\rm (ii)} \ \gamma_{i_{1},j_{1}}\gamma_{i_{2},j_{2}}\gamma_{i_{1},j'_{1}}\cdot c(\gamma,0)\cdot a_{i_{1}}\cdot s_{1,j'_{1}-1,a_{i_{1}}-j'_{1}}\cdot
c(e_{i_{1},1}+e_{i_{1},j'_{1}-1}+e_{i_{1},a_{i_1}-j'_{1}},0)=0}\\
&&\text{if} \ \ j_{1}+j'_{1}> a_{i_{1}}.
\end{eqnarray*}
We shall show that $c(\gamma' +e_{i_{1},1}+e_{i_{2},j_{2}}+e_{i_{1},j_{1}+j'_{1}-1},0)=0$.
One has $\deg(t_{i_{1},j''_{1}})\leq (j_{1}+j'_{1}-1)/a_{i_{1}}$. 
Then, taking the coefficient in front of $t_{i_{1},1}$ in 
$WDVV((i_{1},j''_{1}),(i_{2},j_{2}),(i_{1},j_{1}+j'_{1}-2),(i_{1},1))$, one has
\begin{eqnarray*}
\lefteqn{{\rm (i)} \ (\gamma'_{i_{1},j''_{1}}+1)(\gamma'_{i_{2},j_{2}}+1)(\gamma'_{i_{1},j_{1}+j'_{1}-1}+1)\cdot 
c(\gamma' +e_{i_{1},1}+e_{i_{2},j_{2}}+e_{i_{1},j_{1}+j'_{1}-1},0)\cdot a_{i_{1}}\cdot }\\
&&s_{1,j_{1}+j'_{1}-2,a_{i_{1}}+1-j_{1}-j'_{1}}\cdot c(e_{i_{1},1}+e_{i_{1},j_{1}+j'_{1}-2}+e_{i_{1},a_{i_{1}}+1-j_{1}-j'_{1}},0)\\
&&-\ \ (\gamma''_{i_{1},1}+1)(\gamma''_{i_{2},j_{2}}+1)(\gamma''_{i_{1},a_{i_{1}}-1}+1)\cdot 
c(\gamma' +2e_{i_{1},1}-e_{i_{1},j''_{1}}+e_{i_{2},j_{2}}+e_{i_{1},a_{i_{1}}-1},0)\cdot a_{i_{1}}\cdot \\
&&s_{j_{1}+j'_{1}-1,j''_{1},1}\cdot c(e_{i_{1},j_{1}+j'_{1}-2}+e_{i_{1},j''_{1}}+e_{i_{1},1},0)\\
&&\text{if} \ \ \deg(t_{i_{1},j''_{1}})=\frac{j_{1}+j'_{1}-1}{a_{i_{1}}} \ \text{and} \ \text{where} \ \gamma''=\gamma' +e_{i_{1},1}-e_{i_{1},j''_{1}}=e_{i_1,1},\\
\lefteqn{{\rm (ii)} \ (\gamma'_{i_{1},j''_{1}}+1)(\gamma'_{i_{2},j_{2}}+1)(\gamma'_{i_{1},j_{1}+j'_{1}-1}+1)\cdot 
c(\gamma' +e_{i_{1},1}+e_{i_{2},j_{2}}+e_{i_{1},j_{1}+j'_{1}-1},0)\cdot a_{i_{1}}\cdot }\\
&&s_{1,j_{1}+j'_{1}-2,a_{i_{1}}+1-j_{1}-j'_{1}}\cdot c(e_{i_{1},1}+e_{i_{1},j_{1}+j'_{1}-2}+e_{i_{1},a_{i_{1}}+1-j_{1}-j'_{1}},0)\\
&&\text{if} \ \ \deg(t_{i_{1},j''_{1}})\leq \frac{j_{1}+j'_{1}-2}{a_{i_{1}}}.
\end{eqnarray*}
If one has $c(\gamma' +2e_{i_{1},1}-e_{i_{1},j''_{1}}+e_{i_{2},j_{2}}+e_{i_{1},a_{i_{1}}-1},0 )\neq 0$, one should have 
\[
2\deg(t_{i_{1},1})+\deg(t_{i_{1},a_{i_{1}}-1})+\deg(i_{2},j_{2})\leq 2 \displaystyle\Leftrightarrow  \frac{1}{a_{i_{1}}}+\frac{1}{a_{i_{2}}} \geq 1.
\]
This inequality contradicts the assumption that $a_{i_{1}}\ge 3$ and $a_{i_{2}}\ge 2$.
Then $c(\gamma' +2e_{i_{1},1}-e_{i_{1},j''_{1}}+e_{i_{2},j_{2}}+e_{i_{1},a_{i_{1}}-1},0 )=0$ and 
$c(\gamma' +e_{i_{1},1}+e_{i_{2},j_{2}}+e_{i_{1},j_{1}+j'_{1}-1},0)=0$.
Hence one has $c(\gamma ,0)=0$.
\qed
\end{pf}

\vspace{5pt}
\noindent
\underline{\it Step 1-{\rm (iv)}. The tem
$t^{\gamma}$ has, as factors, $t_{i_1,j_1}$ and $t_{i_{2},j}, t_{i_{2},j'}, t_{i_{2},j''}$.}

\vspace{5pt}
If $a_{i_2}\ge 3$ and some $j\ge 2$, these cases are already dealt with in previous arguments 
in Step 1-{\rm (iii)}. If $a_{i_2}\ge 3$ and $j=j'=j''=1$, one has $\deg(t^{\gamma})>2$.
Therefore, we only have to consider the case $a_{i_2}=2$:
\vspace{-5pt}
\begin{sublem}[Step 1-{\rm (iv)}]\label{sl4-l2-wsep}
Let 
$\gamma \in \ZZ^{\mu_{A}-2}$ be a non-negative element such that $|\gamma|=4$,
$\gamma_{i_{3},j_{3}}=0$ for all $j_{3}\ge 1$ and
$\gamma-e_{i_{1},j_{1}}-2e_{i_{2},1}\ge 0$ for $j_{1}\ge 2$.
If $a_{i_{1}}\ge 3$ and $a_{i_{2}}=2$, then one has $c(\gamma,0)=0$. 
\end{sublem}

\begin{pf}
Taking the coefficient in front of $t^{\gamma-e_{i_{1},j_{1}}-2e_{i_{2},1}}$ 
in $WDVV((i_{2},1),(i_{2},1),(i_{1},j_{1}-1),(i_{1},1))$, 
one has
\[
\gamma_{i_{1},j_{1}}\gamma_{i_{2},1}\cdot (\gamma_{i_{2},1}-1)c(\gamma, 0)\cdot a_{i_{1}}\cdot 
c(e_{i_{1},j_{1}-1}+e_{i_{1},1}+e_{i_{1},a_{i_{1}}-j_{1}},0)=0.
\]
Hence one has $c(\gamma,0)=0$. 
\qed
\end{pf}

\noindent
\vspace{5pt}
\underline{\it Step 2. The term 
$t^{\gamma}$ does not have, as factors, $t_{i,j}$ 
for $j\ge 2$.}

\vspace{10pt}
We shall split Step 2 into the following two cases:

\noindent
Case 1: The term $t^{\gamma}$ has $t_{1,1}, t_{2,1}, t_{3,1}$ as factors.

\noindent
Case 2: The term $t^{\gamma}$ has, as factors, only two parameters $t_{i_1,1}, t_{i_2,1}$
for $i_1\ne i_2$. 
\vspace{-5pt}
\begin{sublem}[Step 2-Case 1]\label{sl2-l1-wsep}
Let $\gamma \in \ZZ^{\mu_{A}-2}$ be a non-negative element such that 
$|\gamma|=4$ and $\gamma=\gamma_{1,1}e_{1,1}+\gamma_{2,1}e_{2,1}+\gamma_{3,1}e_{3,1}$ for $\gamma_{1,1}\gamma_{2,1}\gamma_{3,1}\ne 0$.
Then one has $c(\gamma,0)=0$. 
\end{sublem}

\begin{pf}
By the assumption that $a_{3}\ge 3$,
one has 
\[
\deg (t^{\gamma})\ge 4\frac{a_{l}-1}{a_{l}}\geq 2, 
\]
where $a_{l}=\min \{a_{1},a_{2},a_{3}\}$. 
The first equality is attained if and only if $a_{1}=a_{2}=a_{3}$.
If $a_{1}=a_{2}=a_{3}$, one also has $\deg (t^{\gamma})> 2$.
Hence one has $c(\gamma ,0)=0$.
\qed
\end{pf}

\begin{sublem}[Step 2-Case 2]\label{sl5-l2-wsep}
Let
$\gamma \in \ZZ^{\mu_{A}-2}$ be a non-negative element such that $|\gamma|=4$ and
$\gamma=\gamma_{i_{1},1}e_{i_{1},1}+\gamma_{i_{2},1}e_{i_{2},1}$ for $i_1\ne i_2$ and  $\gamma_{i_{1},1}\gamma_{i_{2},1}\ne 0$.
If $a_{i_{1}}\ge 3$, then one has $c(\gamma,0)=0$. 
\end{sublem}

\begin{pf}
One has 
\[
\deg (t^{\gamma})\ge 4\frac{a_{l}-1}{a_{l}}\geq 2, 
\]
where $a_{l}=\min \{a_{i_{1}},a_{i_{2}}\}$.
The first equality is attained if and only if $a_{1}=a_{2}=a_{3}$.
If $a_{1}=a_{2}=a_{3}$, one also has $\deg (t^{\gamma})> 2$.
Hence one has $c(\gamma,0)=0$.   
\qed
\end{pf}
Therefore we have Lemma~\ref{wsep}.
\qed
\end{pf}

\begin{rem}
Lemma~\ref{wsep} is the weaker version of Proposition~\ref{sep}. Lemma~\ref{wsep}
is, however, enough to show the following Proposition~\ref{lem3.1}, which is necessary
for the proofs of Proposition~\ref{sep}, Proposition~\ref{0reconst} and Proposition~\ref{mreconst}. 
Therefore we, at first, put Lemma~\ref{wsep} as an independent lemma for the readers convenience. 
\end{rem}

In the proof of the following Proposition~\ref{lem3.1}, we shall use the condition {\rm (iv)}
in Theorem~\ref{first} for the case $A=(2,2,2)$ (Sub-Lemma~\ref{s(2,2,2),l2-l3.1}).
For other cases, we can show Proposition~\ref{lem3.1} by other conditions together with 
the weaker condition {\rm (iv')} for the case that $A=(1,2,2)$ and $A=(2,2,a_3), a_3\ge 3$ (Sub-Lemma~\ref{s1,l2-l3.1} and Sub-Lemma~\ref{ishi-lem3.1}) and by other conditions without any further conditions except for the cases that $A=(1,2,2)$ and $A=(2,2,a_3)$, $a_3\ge 3$.
\begin{prop}\label{lem3.1}
If $a_{1}\ge 2$ $($resp. $1=a_{1}<a_{2}\le a_{3}$, $1=a_{1}=a_{2}<a_{3}$$)$, 
$c(\alpha,1)$ with $|\alpha|\le 3$ $($resp. $|\alpha|\le 2$, $|\alpha|\le 1$$)$
is none-zero if and only if $\alpha=e_{1,1}+e_{2,1}+e_{3,1}$ $($resp. $\alpha=e_{2,1}+e_{3,1}$, $\alpha=e_{3,1}$$)$.
In particular, one has $c(e_{1,1}+e_{2,1}+e_{3,1},1)=1$ $($resp. $c(e_{2,1}+e_{3,1},1)=1$, $c(e_{3,1},1)=1$$)$ 
by the condition {\rm (vi)} of Theorem~\ref{first}.
\end{prop}
\begin{pf}

We first show Proposition~\ref{lem3.1} for the case $a_{1}\ge 2$.
We shall split the proof into following two cases.

\begin{lem}[Case 1]\label{lem1-lem3.1}
Let $\gamma \in \ZZ^{\mu_{A}-2}$ be a non-negative element such that $|\gamma|=3$, $\gamma-e_{i,j}\geq 0$ for some $i, j$.
If $a_{i}\ge 3$ and  $j\ge 2$, then one has $c(\gamma,1)=0$.
\end{lem}

\begin{pf}
One has $c(\alpha, 1)=0$ if $|\alpha|\le 2$ because $\deg(t^{\alpha}e^{t_{\mu_{A}}})<2$.
Then, taking the coefficient in front of  
$t^{\gamma-e_{i,j}}e^{t_{\mu_A}}$
in $WDVV((i,1),(i,j-1),\mu_A,\mu_A)$, one has
\[
s_{1,j-1,a_{i}-j}\cdot c(e_{i,1}+e_{i,j-1}+e_{i,a_{i}-j},0)\cdot a_{i} \cdot \gamma_{i,j} \cdot c(\gamma ,1)=0.
\]
Hence one has $c(\gamma ,1)=0$.
\qed
\end{pf}

\begin{lem}[Case 2]\label{lem2-lem3.1}
Let $\gamma \in \ZZ^{\mu_{A}-2}$ be a non-negative element such that
$|\gamma|=3$, $\gamma =\gamma_{1,1}e_{1,1}+\gamma_{2,1}e_{2,1}+\gamma_{3,1}e_{3,1}$ 
and $\gamma_{1,1}\gamma_{2,1}\gamma_{3,1}=0$.
Then one has $c(\gamma ,1)=0$.
\end{lem}

\begin{pf}
We shall split the proof of Lemma \ref{lem2-lem3.1} into following two cases.

\begin{sublem}[Case 2-1]\label{s1,l2-l3.1}
If $A=(2,2,r)$ with $r\ge 3$, 
then one has $c(2e_{1,1}+e_{3,1},1)=c(2e_{2,1}+e_{3,1},1)=0$.
\end{sublem}

\begin{pf}
Note that $c(\alpha,0)=0$ if $|\alpha|=4$, $\alpha-e_{i,1}-e_{3,j_{3}}\ge 0$ for $1\le j_{3}\le a_{3}-1$ and $i=1,2$ by Lemma \ref{wsep}, and $c(\alpha, 1)=0$ if $|\alpha|\le 2$ because $\deg(t^{\alpha}e^{t_{\mu_{A}}})<2$.
Then, taking the coefficient in front of $t_{3,1}e^{t_{\mu_A}}$ in $WDVV((1,1),(1,1),(2,1),(2,1))$, one has
\[
2c(2e_{1,1}+e_{1},0)\cdot 1\cdot 2c(2e_{2,1}+e_{3,1},1)+2c(2e_{2,1}+e_{1},0)\cdot 1\cdot 2c(2e_{1,1}+e_{3,1},1)=0.
\]
One has $c(2e_{1,1}+e_{3,1},1)=c(2e_{2,1}+e_{3,1},1)$ because of the condition {\rm (iv')} in Theorem~\ref{first}.
Hence one has $c(2e_{1,1}+e_{3,1},1)=c(2e_{2,1}+e_{3,1},1)=0$. 
\qed
\end{pf}

\begin{sublem}[Case 2-2]\label{s(2,2,2),l2-l3.1}
Suppose that $A=(2,2,2)$.  Let $\gamma\in \ZZ^{\mu_A-2}$ be a non-negative element such that
$|\gamma|=3$ and $\gamma-2e_{i,1}\ge 0$.
Then one has $c(\gamma,1)=0$.
\end{sublem}

\begin{pf}
Note that $c(\gamma,0)=0$ if $\gamma-e_{i_1,1}-e_{i_2,1}\ge0$ for $i_1\ne i_2$ by the condition {\rm (iv)} in Theorem~\ref{first}.  
Without loss of generality, we can assume that $\gamma_{1,1}=0$.
Taking the coefficient in front of $t_{2,1}^{\gamma_{2,1}}t_{3,1}^{\gamma_{3,1}}e^{t_{\mu_A}}$
in $WDVV((1,1),(1,a_{1}-1),\mu_A,\mu_A)$, one has
\[
c(e_{1}+e_{1,a_{1}-1}+e_{1,1},0)\cdot c(\gamma ,1)=0,
\]
Hence one has $c(\gamma ,1)=0$.
\qed
\end{pf}

\begin{sublem}[Case 2-3]\label{s2,l2-l3.1}
Let $\gamma \in \ZZ^{\mu_{A}-2}$ be a non-negative element such that
$|\gamma|=3$ and $\gamma=\gamma_{1,1}e_{1,1}+\gamma_{2,1}e_{2,1}+\gamma_{3,1}e_{3,1}$ with  
$\gamma_{i,1}=0$ for some $i$.
If $A\ne (2,2,r)$, then one has $c(\gamma,1)=0$. 
\end{sublem}

\begin{pf}
Note that $c(\alpha,0)=0$ if $|\alpha|=4$ and $\alpha-e_{i_{1},j_{1}}-e_{i_{2},j_{2}}\ge 0$ for $i_{1}\neq i_{2}$
by Lemma \ref{wsep}, and $c(\alpha, 1)=0$ if $|\alpha|\le 2$.
Then, taking the coefficient in front of $t_{1,1}^{\gamma_{1,1}}t_{2,1}^{\gamma_{2,1}}t_{3,1}^{\gamma_{3,1}}e^{t_{\mu_A}}$
in $WDVV((i,1),(i,a_{i}-1),\mu_A,\mu_A)$, one has
\[
c(e_{1}+e_{i,a_{i}-1}+e_{i,1},0)\cdot c(\gamma ,1)=0,
\]
Hence one has $c(\gamma ,1)=0$.
\qed
\end{pf}
Then we have Lemma~\ref{lem2-lem3.1}.
\qed
\end{pf}

Next we show Proposition~\ref{lem3.1} for the case  $1=a_{1}<a_{2}\le a_{3}$.
We shall split the proof into following two cases.
\begin{lem}[Case 1]
Let $\gamma \in \ZZ^{\mu_{A}-2}$ be a non-negative element such that $|\gamma|=2$, $\gamma-e_{i,j}\geq 0$ for $2\leq j$.
If $a_{i}\ge 3$, then one has $c(\gamma,1)=0$.
\end{lem}

\begin{pf}
One has $c(\alpha, 1)=0$ if $|\alpha|\le 1$ because $\deg(t^{\alpha}e^{t_{\mu_{A}}})<2$.
Then, taking the coefficient in front of  
$t^{\gamma-e_{i,j}}e^{t_{\mu_A}}$
in $WDVV((i,1),(i,j-1),\mu_A,\mu_A)$, one has
\[
s_{1,j-1,a_{i}-j}\cdot c(e_{i,1}+e_{i,j-1}+e_{i,a_{i}-j},0)\cdot a_{i} \cdot \gamma_{i,j} \cdot c(\gamma ,1)=0.
\]
Hence one has $c(\gamma ,1)=0$.
\qed
\end{pf}

\begin{lem}[Case 2]\label{ishi-lem3.1}
Let $\gamma \in \ZZ^{\mu_{A}-2}$ be a non-negative element such that
$|\gamma|=2$, $\gamma =\gamma_{2,1}e_{2,1}+\gamma_{3,1}e_{3,1}$ 
and $\gamma_{2,1}\gamma_{3,1}=0$.
Then one has $c(\gamma ,1)=0$.
\end{lem}

\begin{pf}
First we shall show Sub-Lemma~\ref{ishi-lem3.1} for the case that $a_{3}\ge 3$.
Note that $c(\alpha,0)=0$ if $|\alpha|=4$ and $\alpha-e_{i_{1},j_{1}}-e_{i_{2},j_{2}}\ge 0$ for $i_{1}\neq i_{2}$
by Lemma \ref{wsep}, and $c(\alpha, 1)=0$ if $|\alpha|\le 1$.
Then, taking the coefficient in front of $t_{1,1}^{\gamma_{1,1}}t_{2,1}^{\gamma_{2,1}}t_{3,1}^{\gamma_{3,1}}e^{t_{\mu_A}}$
in $WDVV((i,1),(i,a_{i}-1),\mu_A,\mu_A)$, one has
\[
c(e_{1}+e_{i,a_{i}-1}+e_{i,1},0)\cdot c(\gamma ,1)=0,
\]
Hence one has $c(\gamma ,1)=0$.

Next we shall show this claim for the case that $a_{3}=2$, i.e., $A=(1,2,2)$.
Taking the coefficient in front of $e^{t_{\mu_A}}$ in $WDVV((2,1),(2,1),(3,1),(3,1))$, one has
\[
2c(2e_{2,1}+e_{1},0)\cdot 1\cdot 2c(2e_{3,1},1)+2c(2e_{3,1}+e_{1},0)\cdot 1\cdot 2c(2e_{2,1},1)=0.
\]
One has $c(2e_{3,1},1)=c(2e_{2,1},1)$ because of the condition {\rm (iv')} in Theorem~\ref{first}.
Hence one has $c(2e_{2,1},1)=c(2e_{3,1},1)=0$. 
\qed
\end{pf}

Finally we show Proposition~\ref{lem3.1} for the case  $1=a_{1}=a_{2}<a_{3}$.

\begin{lem}\label{ishi}
If $\gamma=e_{3,j}\geq 0$ for $2\leq j$,  $c(\gamma ,1)=0$.
\end{lem}

\begin{pf}
One has $\deg(t^{\gamma}e^{t_{\mu_{A}}})<2$. Therefore Lemma~\ref{ishi} holds.
\qed
\end{pf}
Therefore we have Proposition~\ref{lem3.1}.
\qed
\end{pf}

The following Proposition~\ref{sep} is nothing but the condition {(iv)} in Theorem~\ref{first}.
However, for the later convenience, 
we shall show Proposition~\ref{sep} is derived from the other conditions together with the weaker condition {\rm (iv)} for the
cases that $A=(1,2,2)$ and $A=(2,2,a_3), a_3\ge 3$ (Sub-Lemma~\ref{sl6-l2-sep} and Sub-Lemma~\ref{sl7-l2-sep}) 
and by other conditions without any further conditions except for the cases that $A=(1,2,2)$ and $A=(2,2,a_3)$, $a_3\ge 3$.
\begin{prop}\label{sep}
Suppose that $A\ne (2,2,2)$.
For a non-negative $\beta\in \ZZ^{\mu_A-2}$, one has
\[
c\left(\beta+\sum_{k=1}^3e_{i_k,j_{k}},0\right)\ne 0
\]
only if $i_1=i_2=i_3$.
\end{prop}

\begin{pf}
We will prove Proposition~\ref{sep} by the induction on the length.
By Proposition \ref{lem3}, one has $c(\alpha,0)=0$ if $|\alpha|=3$, $\alpha-e_{i_{1},j_{1}}-e_{i_{2},j_{2}}\geq 0$ for $i_{1}\neq i_{2}$.
Assume that $c(\alpha,0)=0$ if $|\alpha|\le k+3$, $\alpha-e_{i_{1},j_{1}}-e_{i_{2},j_{2}}\geq 0$ for $i_{1}\neq i_{2}$.
Under this assumption,
we will prove that $c(\alpha,0)=0$ if $|\alpha|=k+4$, 
$\alpha-e_{i_{1},j_{1}}-e_{i_{2},j_{2}}\geq 0$ for $i_{1}\neq i_{2}$. 

We shall split the proof into the following three steps:

\vspace{5pt}
\noindent
Step 1: We shall consider the case that the term $t^{\gamma}$ has, as a factor,
$t_{i_1,j}$ such that $a_{i_1}\ge 3$ and $j\ge 2$ for some $j$.
We shall split, moreover, Step 1 into following four cases: 
\begin{enumerate}
\item The term $t^{\gamma}$ has $t_{1,j_1}, t_{2,j_2}, t_{3,j_3}$ as factors.
\item The term $t^{\gamma}$ has $t_{i,j}$ and $t_{i,j'}$ as factors for each $i=i_{1}, i_{2}$.
\item The term $t^{\gamma}$ has $t_{i_1,j_1}, t_{i_1,j'_1}, t_{i_1,j''_1}$ and only $t_{i_2,j_2}$ as factors.
\item The term $t^{\gamma}$ has, as factors, $t_{i_1,j_1}$ and $t_{i_{2},j}, t_{i_{2},j'}, t_{i_{2},j''}$.
\end{enumerate}

\vspace{5pt}
\noindent
Step 2: We shall consider the case that $a_{i_1}\ge 3$ and $t^{\gamma}$ does not have, as factors, $t_{i,j}$ for $j\ge 2$. 

\vspace{5pt}
\noindent
Step 3: We shall consider the case that $a_{1}=a_{2}=2$ and $a_3\ge 3$.

\vspace{5pt}
\noindent
Step 4: We shall consider the case that $a_{1}=1, a_{2}=a_3=2$.

\vspace{5pt}
\noindent
\underline{\it Step 1-{\rm (i)}. The term $t^{\gamma}$ has $t_{1,j_1}, t_{2,j_2}, t_{3,j_3}$ as factors. }

\begin{sublem}[Step 1-(i)]\label{sl1-l1-sep}
Assume that $a_{i_{3}}\ge 3$ and $c(\alpha,0)=0$ if $|\alpha|\le k+3$ and $\alpha-e_{i_{1},j_{1}}-e_{i_{2},j_{2}}\geq 0$ for $i_{1}\neq i_{2}$.
Let $\gamma \in \ZZ^{\mu_{A}-2}$ be a non-negative element such that 
$|\gamma|=k+4$, $\gamma -e_{i_{1},j_{1}}-e_{i_{2},j_{2}}-e_{i_{3},j_{3}}\ge 0$ for distinct $i_{1}$, $i_{2}$, $i_{3}$ and $j_{3}\ge 2$.
Then one has $c(\gamma,0)=0$. 
\end{sublem}

\begin{pf}
Taking the coefficient in front of $t^{\gamma -e_{i_{1},j_{1}}-e_{i_{2},j_{2}}-e_{i_{3},j_{3}}}$ in $WDVV((i_{1},j_{1}),(i_{2},j_{2}),(i_{3},j_{3}-1),(i_{3},1))$,
one has
\[
\gamma_{i_{1},j_{1}}\gamma_{i_{2},j_{2}}\gamma_{i_{3},j_{3}}\cdot c(\gamma,0)\cdot a_{i_{3}}\cdot 
s_{a_{i_{3}}-j_{3},j_{3}-1,1}\cdot c(e_{i_{3},a_{i_{3}}-j_{3}}+e_{i_{3},j_{3}-1}+e_{i_{3},1},0)=0.
\]
Hence one has $c(\gamma ,0)=0$.
\qed
\end{pf}

\vspace{5pt}
\noindent
\underline{\it Step 1-{\rm (ii)}.  The term
$t^{\gamma}$ has $t_{i,j}$ and $t_{i,j'}$ as factors for each $i=i_{1}, i_{2}$.}

\begin{sublem}[Step 1-(ii)]\label{sl1-l2-sep}
Assume that $a_{i_{1}}\geq 3$ and $c(\alpha,0)=0$ if $|\alpha|\le k+3$ and $\alpha-e_{i_{1},j_{1}}-e_{i_{2},j_{2}}\geq 0$ for $i_{1}\neq i_{2}$. 
Let $\gamma \in \ZZ^{\mu_{A}-2}$ be a non-negative element such that $|\gamma|=k+4$,
$\gamma_{i_{3},j_{3}}=0$ for $1\leq \forall j_{3}\leq a_{i_{3}}-1$ and
$\gamma-e_{i_{1},j_{1}}-e_{i_{1},j'_{1}}-e_{i_{2},j_{2}}-e_{i_{2},j'_{2}}\geq 0$ for $j'_{1}\ge 2$ and some $j'_{2}$. 
Then one has $c(\gamma,0)=0$.
\end{sublem}

\begin{pf}
Taking the coefficient in front of $t^{\gamma-e_{i_{1},j_{1}}-e_{i_{2},j_{2}}-e_{i_{1},j'_{1}}}$ in $WDVV((i_{2},j_{2}),(i_{1},j_{1}),(i_{1},1),(i_{1},j'_{1}-1))$, 
one has
\begin{eqnarray*}
\lefteqn{{\rm (i)} \ \gamma_{i_{1},j_{1}}\gamma_{i_{2},j_{2}}\gamma_{i_{1},j'_{1}}\cdot c(\gamma,0)\cdot a_{i_{1}}\cdot 
s_{1,j'_{1}-1,a_{i_{1}}-j'_{1}}\cdot c(e_{i_{1},1}+e_{i_{1},j'_{1}-1}+e_{i_{1},a_{i_{1}}-j'_{1}},0)}\\
&&-(\gamma'_{i_{1},1}+1)(\gamma'_{i_{2},j_{2}}+1)(\gamma'_{i_{1},j_{1}+j'_{1}-1}+1)\cdot c(\gamma' +e_{i_{1},1}+e_{i_{2},j_{2}}+e_{i_{1},j_{1}+j'_{1}-1},0)
\cdot a_{i_{1}}\cdot \\
&&s_{j_1,j'_{1}-1,a_{i_{1}}+1-j_{1}-j'_{1}}\cdot c(e_{i_{1},j_{1}}+e_{i_{1},j'_{1}-1}+e_{i_{1},a_{i_{1}}+1-j_{1}-j'_{1}},0)=0\\ 
&&\text{if} \ \ 3\leq j_{1}+j'_{1}\leq a_{i_{1}} \ \text{and} \ \text{where} \ \gamma'=\gamma-e_{i_{1},j_{1}}-e_{i_{2},j_{2}}-e_{i_{1},j'_{1}},\\ 
\lefteqn{{\rm (ii)} \ \gamma_{i_{1},j_{1}}\gamma_{i_{2},j_{2}}\gamma_{i_{1},j'_{1}}\cdot c(\gamma,0)\cdot a_{i_{1}}\cdot s_{1,j'_{1}-1,a_{i_{1}}-j'_{1}}\cdot
c(e_{i_{1},1}+e_{i_{1},j'_{1}-1}+e_{i_{1},a_{i_1}-j'_{1}},0)=0}\\
&&\text{if} \ \ j_{1}+j'_{1}> a_{i_{1}}.
\end{eqnarray*}
We shall show that $c(\gamma' +e_{i_{1},1}+e_{i_{2},j_{2}}+e_{i_{1},j_{1}+j'_{1}-1},0)=0$.
One has $\deg(t^{\gamma'}t_{i_{2},j_{2}})=(j_{1}+j'_{1})/a_{i_{1}}$, i.e, $\deg(t_{i_{2},j_{2}}t_{i_{2},j'_{2}})\leq (j_{1}+j'_{1})/a_{i_{1}}\le 1$. 
Then, taking the coefficient in front of $t^{\gamma' +e_{i_{1},1}-e_{i_{2},j'_{2}}}$ in $WDVV((i_{2},j_{2}),(i_{2},j'_{2}),(i_{1},j_{1}+j'_{1}-2),(i_{1},1))$, 
one has
\begin{multline*}
(\gamma'_{i_{2},j_{2}}+1)(\gamma'_{i_{2},j'_{2}}+1)(\gamma'_{i_{1},j_{1}+j'_{1}-1}+1)\cdot c(\gamma' +e_{i_{1},1}+e_{i_{2},j_{2}}+e_{i_{1},j_{1}+j'_{1}-1},0)
\cdot a_{i_{1}} \cdot \\
s_{1,j_{1}+j'_{1}-2,a_{i_{1}}+1-j_{1}-j'_{1}}\cdot c(e_{i_{1},1}+e_{i_{1},j_{1}+j'_{1}-2}+e_{i_{1},a_{i_{1}}+1-j_{1}-j'_{1}},0)=0.
\end{multline*}
By Proposition \ref{lem3}, $c(\gamma' +e_{i_{1},1}+e_{i_{2},j_{2}}+e_{i_{1},j_{1}+j'_{1}-1},0)=0$.
Hence one has $c(\gamma,0)=0$.
\qed
\end{pf}

\vspace{5pt}
\noindent
\underline{\it Step 1-{\rm (iii)}. The term
$t^{\gamma}$ has $t_{i_1,j_1}, t_{i_1,j'_1}, t_{i_1,j''_1}$ and $t_{i_2,j_2}$ as factors.}

\vspace{10pt}
We shall split Step 1-{\rm (iii)} into Case 1 ($j_2\ge 2$) and Case 2 ($j_2=1$).
\vspace{-5pt}
\begin{sublem}[Step 1-(iii)-Case 1]\label{sl2-l2-sep}
Assume that $a_{i_{1}}\ge 3$, $a_{i_{2}}\ge 3$ and $c(\alpha,0)=0$ 
if $|\alpha|\le k+3$ and $\alpha-e_{i_{1},j_{1}}-e_{i_{2},j_{2}}\geq 0$ for $i_{1}\neq i_{2}$. 
Let $\gamma \in \ZZ^{\mu_{A}-2}$ be a non-negative element such that $|\gamma|=k+4$,
$\gamma_{i_{3},j_{3}}=0$ for $1\leq \forall j_{3}\leq a_{i_{3}}-1$ and
$\gamma-e_{i_{1},j_{1}}-e_{i_{1},j'_{1}}-e_{i_{1},j''_{1}}-e_{i_{2},j_{2}}\geq 0$ for $j'_{1}\ge 2$ and $j_{2}\ge 2$. 
Then one has $c(\gamma,0)=0$.
\end{sublem}

\begin{pf}
Taking the coefficient in front of $t^{\gamma-e_{i_{1},j_{1}}-e_{i_{2},j_{2}}-e_{i_{1},j'_{1}}}$ in $WDVV((i_{2},j_{2}),(i_{1},j_{1}),(i_{1},1),(i_{1},j'_{1}-1))$, 
one has
\begin{eqnarray*}
\lefteqn{{\rm (i)} \ \gamma_{i_{1},j_{1}}\gamma_{i_{2},j_{2}}\gamma_{i_{1},j'_{1}}\cdot c(\gamma,0)\cdot a_{i_{1}}\cdot 
s_{1,j'_{1}-1,a_{i_{1}}-j'_{1}}\cdot c(e_{i_{1},1}+e_{i_{1},j'_{1}-1}+e_{i_{1},a_{i_{1}}-j'_{1}},0)}\\
&&-(\gamma'_{i_{1},1}+1)(\gamma'_{i_{2},j_{2}}+1)(\gamma'_{i_{1},j_{1}+j'_{1}-1}+1)\cdot c(\gamma' +e_{i_{1},1}+e_{i_{2},j_{2}}+e_{i_{1},j_{1}+j'_{1}-1},0)
\cdot a_{i_{1}}\cdot \\
&&s_{j_1,j'_{1}-1,a_{i_{1}}+1-j_{1}-j'_{1}}\cdot c(e_{i_{1},j_{1}}+e_{i_{1},j'_{1}-1}+e_{i_{1},a_{i_{1}}+1-j_{1}-j'_{1}},0)=0\\ 
&&\text{if} \ \ 3\leq j_{1}+j'_{1}\leq a_{i_{1}} \ \text{and} \ \text{where} \ \gamma'=\gamma-e_{i_{1},j_{1}}-e_{i_{2},j_{2}}-e_{i_{1},j'_{1}},\\ 
\lefteqn{{\rm (ii)} \ \gamma_{i_{1},j_{1}}\gamma_{i_{2},j_{2}}\gamma_{i_{1},j'_{1}}\cdot c(\gamma,0)\cdot a_{i_{1}}\cdot s_{1,j'_{1}-1,a_{i_{1}}-j'_{1}}\cdot
c(e_{i_{1},1}+e_{i_{1},j'_{1}-1}+e_{i_{1},a_{i_1}-j'_{1}},0)=0}\\
&&\text{if} \ \ j_{1}+j'_{1}> a_{i_{1}}.
\end{eqnarray*}
We shall show that $c(\gamma' +e_{i_{1},1}+e_{i_{2},j_{2}}+e_{i_{1},j_{1}+j'_{1}-1},0)=0$.
One has $\deg(t_{i_{1},j''_{1}})\leq (j_{1}+j'_{1}-1)/a_{i_{1}}$. Then, taking 
the coefficient in front of $t^{\gamma' +e_{i_{1},1}-e_{i_{1},j''_{1}}}$ in $WDVV((i_{1},j''_{1}),(i_{1},j_{1}+j'_{1}-1),(i_{2},j_{2}-1),(i_{2},1))$,
one has 
\begin{multline*}
\lefteqn{(\gamma'_{i_{1},j''_{1}}+1)(\gamma'_{i_{1},j_{1}+j'_{1}-1}+1)(\gamma'_{i_{2},j_{2}}+1)\cdot c(\gamma 
+e_{i_{1},1}+e_{i_{2},j_{2}}+e_{i_{1},j_{1}+j'_{1}-1},0)
\cdot a_{i_{2}}\cdot}\\
s_{1,j_{2}-1,a_{i_{2}}-j_{2}}\cdot c(e_{i_{2},1}+e_{i_{2},j_{2}-1}+e_{i_{2},a_{i_{2}}-j_{2}},0)=0.
\end{multline*}
By Proposition \ref{lem3}, $c(\gamma' +e_{i_{1},1}+e_{i_{2},j_{2}}+e_{i_{1},j_{1}+j'_{1}-1},0)=0$.
Hence one has $c(\gamma,0)=0$.
\qed
\end{pf}

\begin{sublem}[Step 1-(iii)-Case 2]\label{sl3-l2-sep}
Assume that $a_{i_{1}}\ge 3$ and $c(\alpha,0)=0$ 
if $|\alpha|\le k+3$ and $\alpha-e_{i_{1},j_{1}}-e_{i_{2},j_{2}}\geq 0$ for $i_{1}\neq i_{2}$. 
Let $\gamma \in \ZZ^{\mu_{A}-2}$ be a non-negative element such that $|\gamma|=k+4$,
$\gamma_{i_{3},j_{3}}=0$ for $1\leq \forall j_{3}$ and
$\gamma-e_{i_{1},j_{1}}-e_{i_{1},j'_{1}}-e_{i_{1},j''_{1}}-e_{i_{2},j_{2}}\geq 0$ for $j'_{1}\ge 2$ and $j_{2}=1$. 
Then one has $c(\gamma,0)=0$.
\end{sublem}

\begin{pf}
Taking the coefficient in front of $t^{\gamma-e_{i_{1},j_{1}}-e_{i_{2},j_{2}}-e_{i_{1},j'_{1}}}$ in $WDVV((i_{2},j_{2}),(i_{1},j_{1}),(i_{1},1),(i_{1},j'_{1}-1))$, 
one has
\begin{eqnarray*}
\lefteqn{{\rm (i)} \ \gamma_{i_{1},j_{1}}\gamma_{i_{2},j_{2}}\gamma_{i_{1},j'_{1}}\cdot c(\gamma,0)\cdot a_{i_{1}}\cdot 
s_{1,j'_{1}-1,a_{i_{1}}-j'_{1}}\cdot c(e_{i_{1},1}+e_{i_{1},j'_{1}-1}+e_{i_{1},a_{i_{1}}-j'_{1}},0)}\\
&&-(\gamma'_{i_{1},1}+1)(\gamma'_{i_{2},j_{2}}+1)(\gamma'_{i_{1},j_{1}+j'_{1}-1}+1)\cdot c(\gamma' +e_{i_{1},1}+e_{i_{2},j_{2}}+e_{i_{1},j_{1}+j'_{1}-1},0)
\cdot a_{i_{1}}\cdot \\
&&s_{j_1,j'_{1}-1,a_{i_{1}}+1-j_{1}-j'_{1}}\cdot c(e_{i_{1},j_{1}}+e_{i_{1},j'_{1}-1}+e_{i_{1},a_{i_{1}}+1-j_{1}-j'_{1}},0)=0\\ 
&&\text{if} \ \ 3\leq j_{1}+j'_{1}\leq a_{i_{1}} \ \text{and} \ \text{where} \ \gamma'=\gamma-e_{i_{1},j_{1}}-e_{i_{2},j_{2}}-e_{i_{1},j'_{1}},\\ 
\lefteqn{{\rm (ii)} \ \gamma_{i_{1},j_{1}}\gamma_{i_{2},j_{2}}\gamma_{i_{1},j'_{1}}\cdot c(\gamma,0)\cdot a_{i_{1}}\cdot s_{1,j'_{1}-1,a_{i_{1}}-j'_{1}}\cdot
c(e_{i_{1},1}+e_{i_{1},j'_{1}-1}+e_{i_{1},a_{i_{1}}-j'_{1}},0)=0}\\
&&\text{if} \ \ j_{1}+j'_{1}> a_{i_{1}}.
\end{eqnarray*}
We shall show that $c(\gamma' +e_{i_{1},1}+e_{i_{2},j_{2}}+e_{i_{1},j_{1}+j'_{1}-1},0)=0$.
One has $\deg(t_{i_{1},j''_{1}})\leq (j_{1}+j'_{1}-1)/a_{i_{1}}$. Then,
taking the coefficient in front of $t^{\gamma' +e_{i_{1},1}-e_{i_{1},j''_{1}}}$ in 
$WDVV((i_{1},j''_{1}),(i_{2},j_{2}),(i_{1},j_{1}+j'_{1}-2),(i_{1},1))$, one has
\begin{eqnarray*}
\lefteqn{{\rm (i)} \ (\gamma'_{i_{1},j''_{1}}+1)(\gamma'_{i_{2},j_{2}}+1)(\gamma'_{i_{1},j_{1}+j'_{1}-1}+1)\cdot 
c(\gamma' +e_{i_{1},1}+e_{i_{2},j_{2}}+e_{i_{1},j_{1}+j'_{1}-1},0)\cdot a_{i_{1}}\cdot }\\
&&s_{1,j_{1}+j'_{1}-2,a_{i_{1}}+1-j_{1}-j'_{1}}\cdot c(e_{i_{1},1}+e_{i_{1},j_{1}+j'_{1}-2}+e_{i_{1},a_{i_{1}}+1-j_{1}-j'_{1}},0)\\
&&-\ \ (\gamma''_{i_{1},1}+1)(\gamma''_{i_{2},j_{2}}+1)(\gamma''_{i_{1},a_{i_{1}}-1}+1)\cdot 
c(\gamma' +2e_{i_{1},1}-e_{i_{1},j''_{1}}+e_{i_{2},j_{2}}+e_{i_{1},a_{i_{1}}-1},0)\cdot a_{i_{1}}\cdot \\
&&s_{j_{1}+j'_{1}-1,j''_{1},1}\cdot c(e_{i_{1},j_{1}+j'_{1}-2}+e_{i_{1},j''_{1}}+e_{i_{1},1},0)=0\\
&&\text{if} \ \ \deg(t_{i_{1},j''_{1}})=\frac{j_{1}+j'_{1}-1}{a_{i_{1}}} \ \text{and} \ \text{where} \ \gamma''=\gamma' +e_{i_{1},1}-e_{i_{1},j''_{1}},\\
\lefteqn{{\rm (ii)} \ (\gamma'_{i_{1},j''_{1}}+1)(\gamma'_{i_{2},j_{2}}+1)(\gamma'_{i_{1},j_{1}+j'_{1}-1}+1)\cdot 
c(\gamma' +e_{i_{1},1}+e_{i_{2},j_{2}}+e_{i_{1},j_{1}+j'_{1}-1},0)\cdot a_{i_{1}}\cdot }\\
&&s_{1,j_{1}+j'_{1}-2,a_{i_{1}}+1-j_{1}-j'_{1}}\cdot c(e_{i_{1},1}+e_{i_{1},j_{1}+j'_{1}-2}+e_{i_{1},a_{i_{1}}+1-j_{1}-j'_{1}},0)=0\\
&&\text{if} \ \ \deg(t_{i_{1},j''_{1}})\leq \frac{j_{1}+j'_{1}-2}{a_{i_{1}}}.
\end{eqnarray*}
If $c(\gamma' +2e_{i_{1},1}-e_{i_{1},j''_{1}}+e_{i_{2},j_{2}}+e_{i_{1},a_{i_{1}}-1},0 )\neq 0$, 
one should have 
\[
2\deg(t_{i_{1},1})+\deg(t_{i_{1},a_{i_{1}}-1})+\deg(i_{2},j_{2})\leq 2 \Leftrightarrow  \frac{1}{a_{i_{1}}}+\frac{1}{a_{i_{2}}} \geq 1.
\]
This inequality contradicts the assumption that $a_{i_{1}}\ge 3$ and $a_{i_{2}}=2$.
Then $c(\gamma' +2e_{i_{1},1}-e_{i_{1},j''_{1}}+e_{i_{2},j_{2}}+e_{i_{1},a_{i_{1}}-1},0 )=0$ and 
$c(\gamma' +e_{i_{1},1}+e_{i_{2},j_{2}}+e_{i_{1},j_{1}+j'_{1}-1},0)=0$.
Hence one has $c(\gamma ,0)=0$.
\qed
\end{pf}

\vspace{5pt}
\noindent
\underline{\it Step 1-{\rm (iv)}. The tem
$t^{\gamma}$ has, as factors, $t_{i_1,j_1}$ and $t_{i_{2},j}, t_{i_{2},j'}, t_{i_{2},j''}$.}

\vspace{5pt}
If $a_{i_2}\ge 3$ and some $j\ge 2$, these cases are already dealt with in previous arguments 
in Step 1-{\rm (iii)}. If $a_{i_2}\ge 3$ and $j=j'=j''=1$, one has $\deg(t^{\gamma})>2$.
Therefore, we only have to consider the case $a_{i_2}=2$:
\vspace{-5pt}
\begin{sublem}[Step 1-(iv)]\label{sl4-l2-sep}
Assume that $a_{i_{1}}\ge 3$, $a_{i_{2}}=2$ and $c(\alpha,0)=0$ 
if $|\alpha|\le k+3$ and $\alpha-e_{i_{1},j_{1}}-e_{i_{2},j_{2}}\geq 0$ for $i_{1}\neq i_{2}$. 
Let $\gamma \in \ZZ^{\mu_{A}-2}$ be a non-negative element such that $|\gamma|=k+4$,
$\gamma_{i_{3},j_{3}}=0$ for all $j_{3}\ge 1$ and
$\gamma-e_{i_{1},j_{1}}-2e_{i_{2},1}\ge 0$ for $j_{1}\ge 2$.
Then one has $c(\gamma,0)=0$. 
\end{sublem}

\begin{pf}
Taking the coefficient in front of $t^{\gamma-e_{i_{1},j_{1}}-2e_{i_{2},1}}$ 
in $WDVV((i_{2},1),(i_{2},1),(i_{1},j_{1}-1),(i_{1},1))$, 
one has
\[
\gamma_{i_{1},j_{1}}\gamma_{i_{2},1}\cdot (\gamma_{i_{2},1}-1)c(\gamma, 0)\cdot a_{i_{1}}\cdot 
c(e_{i_{1},j_{1}-1}+e_{i_{1},1}+e_{i_{1},a_{i_{1}}-j_{1}},0)=0.
\]
Hence one has $c(\gamma,0)=0$. 
\qed
\end{pf}

\vspace{5pt}
\noindent
\underline{\it Step 2. $a_{i_1}\ge 3$ and $t^{\gamma}$ does not have, as factors, $t_{i,j}$ 
for $j\ge 2$.}

\vspace{10pt}
We shall split Step 2 into the following two cases:

\noindent
Case 1: The term $t^{\gamma}$ has $t_{1,1}, t_{2,1}, t_{3,1}$ as factors.

\noindent
Case 2: The term $t^{\gamma}$ has, as factors, only two parameters $t_{i_1,1}, t_{i_2,1}$
for $i_1\ne i_2$. 
\vspace{-5pt}
\begin{sublem}[Step 2-Case 1]\label{sl2-l1-sep}
Assume that $a_{i_{3}}\ge 3$ and $c(\alpha,0)=0$ if $|\alpha|\le k+3$ and $\alpha-e_{i_{1},j_{1}}-e_{i_{2},j_{2}}\geq 0$ for $i_{1}\neq i_{2}$.
Let $\gamma \in \ZZ^{\mu_{A}-2}$ be a non-negative element such that 
$|\gamma|=k+4$ and 
$\gamma=\gamma_{1,1}e_{1,1}+\gamma_{2,1}e_{2,1}+\gamma_{3,1}e_{3,1}$ for 
$\gamma_{1,1}\gamma_{2,1}\gamma_{3,1}\ne 0$.
Then one has $c(\gamma,0)=0$. 
\end{sublem}

\begin{pf}
One has 
\[
\deg (t^{\gamma})\ge 4\frac{a_{l}-1}{a_{l}}\geq 2, 
\]
where $a_{l}=\min \{a_{1},a_{2},a_{3}\}$. 
The first equality is attained if and only if $a_{1}=a_{2}=a_{3}$.
If $a_{1}=a_{2}=a_{3}$, one also has $\deg (t^{\gamma})> 2$.
Hence one has $c(\gamma ,0)=0$.
\qed
\end{pf}

\begin{sublem}[Step 2-Case 2]\label{sl5-l2-sep}
Assume that $a_{i_{1}}\ge 3$ and $c(\alpha,0)=0$ 
if $|\alpha|\le k+3$ and $\alpha-e_{i_{1},j_{1}}-e_{i_{2},j_{2}}\geq 0$ for $i_{1}\neq i_{2}$. 
Let $\gamma \in \ZZ^{\mu_{A}-2}$ be a non-negative element such that $|\gamma|=k+4$
and $\gamma=\gamma_{i_{1},1}e_{i_{1},1}+\gamma_{i_{2},1}e_{i_{2},1}$.
Then one has $c(\gamma,0)=0$. 
\end{sublem}

\begin{pf}
One has 
\[
\deg (t^{\gamma})\ge 4\frac{a_{l}-1}{a_{l}}\geq 2, 
\]
where $a_{l}=\min \{a_{i_{1}},a_{i_{2}}\}$.
The first equality is attained if and only if $a_{1}=a_{2}=a_{3}$.
If $a_{1}=a_{2}=a_{3}$, one also has $\deg (t^{\gamma})> 2$.
Hence one has $c(\gamma,0)=0$.   
\qed
\end{pf}

\vspace{5pt}
\noindent
\underline {\it Step 3. $a_{1}=a_{2}=2$ and $a_3\ge 3$.}

\vspace{5pt}
If $t^{\gamma}$ has, as a factor, $t_{3,j}$ for $j\ge 1$, we have $c(\gamma,0)=0$ by previous arguments in Step 1 and Step 2. Therefore we only have to show the following Sub-Lemma~\ref{sl6-l2-sep}.
\vspace{-5pt}
\begin{sublem}[Step 3]\label{sl6-l2-sep}
Let $A$ be the type $(2,2,r)$.
If $r\ge 3$, one has 
\[
\begin{cases}
{\rm (i)} \ \ c(2e_{1,1}+2e_{2,1},0)=0,\\
{\rm (ii)} \ \ c(3e_{1,1}+e_{2,1},0)=c(3e_{2,1}+e_{1,1},0)=0.
\end{cases}
\]
\end{sublem}

\begin{pf}
Note that $c(\gamma,1)\neq 0$ with $|\gamma|=3$
if and only if $\gamma=e_{1,1}+e_{2,1}+e_{3,1}$ by Proposition~\ref{lem3.1}, and
$c(\gamma,0)=0$ if $\gamma \in \ZZ^{\mu_{A}-2}$ be a non-negative element such that 
$|\gamma|=4$ and $\gamma -e_{i_{1},j_{1}}-e_{3,j_{3}}\ge 0$ for $i_{1}\neq 3$
by Lemma~\ref{wsep}.

First, we shall show that $c(2e_{1,1}+2e_{2,1},0)=0$. 
Taking the coefficient in front of $t_{2,1}^{2}t_{3,1}e^{t_{\mu_A}}$ in $WDVV((1,1),(2,1),\mu_A,\mu_A)$,
one has
\[
4\cdot c(2e_{1,1}+2e_{2,1},0)\cdot 2\cdot 1=0. 
\]
Hence one has $c(2e_{1,1}+2e_{2,1},0)=0$.

Next, we shall show that $c(3e_{1,1}+e_{2,1},0)=c(3e_{2,1}+e_{1,1},0)=0$.
Taking the coefficient in front of $t^{2}_{1,1}t_{3,1}e^{t_{\mu_A}}$ in $WDVV((1,1),(1,1),\mu_A,\mu_A)$, 
one has
\[
6\cdot c(3e_{1,1}+e_{2,1},0)\cdot a_{2}\cdot c(e_{1,1}+e_{2,1}+e_{3,1},1)=0. 
\]
Thus one has $c(3e_{1,1}+e_{2,1},0)=0$.
The same argument shows $c(3e_{2,1}+e_{1,1},0)=0$. 
\qed
\end{pf}

\vspace{5pt}
\noindent
\underline {\it Step 4. $a_{1}=1, a_{2}=a_3=2$.}

\vspace{5pt}
\begin{sublem}[Step 4]\label{sl7-l2-sep}
Let $A$ be the type $(1,2,2)$.
Then one has 
\[
\begin{cases}
{\rm (i)} \ \ c(2e_{2,1}+2e_{3,1},0)=0,\\
{\rm (ii)} \ \ c(3e_{2,1}+e_{3,1},0)=c(3e_{3,1}+e_{2,1},0)=0.
\end{cases}
\]
\end{sublem}

\begin{pf}
Note that $c(\gamma,1)\neq 0$ with $|\gamma|=2$
if and only if $\gamma=e_{2,1}+e_{3,1}$ by Proposition~\ref{lem3.1}.

First, we shall show that $c(2e_{2,1}+2e_{3,1},0)=0$. 
Taking the coefficient in front of $t_{2,1}^{2}e^{t_{\mu_A}}$ in $WDVV((2,1),(3,1),\mu_A,\mu_A)$,
one has
\[
4\cdot c(2e_{2,1}+2e_{3,1},0)\cdot a_3 \cdot c(e_{2,1}+e_{3,1},1)=0. 
\]
Hence one has $c(2e_{2,1}+2e_{3,1},0)=0$.

Next, we shall show that $c(3e_{2,1}+e_{3,1},0)=c(3e_{3,1}+e_{2,1},0)=0$.
Taking the coefficient in front of $t_{3,1}^{2}e^{t_{\mu_A}}$ in $WDVV((2,1),(2,1),\mu_A,\mu_A)$, 
one has
\[
6\cdot c(3e_{2,1}+e_{3,1},0)\cdot a_{2}\cdot c(e_{2,1}+e_{3,1},1)=0. 
\]
Thus one has $c(3e_{2,1}+e_{3,1},0)=0$.
The same argument shows $c(3e_{3,1}+e_{2,1},0)=0$. 
\qed
\end{pf}

Therefore we have Proposition~\ref{sep}.
\qed
\end{pf}

\begin{cor}\label{lem4-1}
If $a_{i_{1}}\ge 3$, then one has 
\[
c(2e_{i_{1},1}+2e_{i_{1},a_{i_{1}}-1},0)=-\frac{1}{4a_{i_{1}}^{2}}. 
\]
\end{cor}

\begin{pf}
By Proposition \ref{sep},
one has $c(e_{i_{1},1}+e_{i_{1},a_{i_{1}}-1}+e_{i_{2},j_{2}}+e_{i_{3},j_{3}}, 0)=0$ if $i_{1}\neq i_{2}$. Then,
taking the coefficient in front of $t_{1,1}t_{2,1}t_{3,1}e^{t_{\mu_A}}$ in $WDVV((i_{1},1),(i_{1},a_{i_{1}}-1),\mu_A ,\mu_A)$, 
one has
\begin{multline*}
c(e_{1}+e_{i_{1},1}+e_{i_{1},a_{i_{1}}-1},0)\cdot 1\cdot c(e_{1,1}+e_{2,1}+e_{3,1},1)+\\
4\cdot c(2e_{i_{1},1}+2e_{i_{1},a_{i_{1}}-1},0)\cdot a_{i_{1}}\cdot c(e_{1,1}+e_{2,1}+e_{3,1},1)=0.
\end{multline*}
By the conditions {\rm (ii)} and {\rm (vi)} in Theorem~\ref{first},
$c(e_{1}+e_{i_{1},1}+e_{i_{1},a_{i_{1}}-1},0)=1/a_{i_{1}}$ and
$c(e_{1,1}+e_{2,1}+e_{3,1},1)=1$.
Hence one has $c(2e_{i_{1},1}+2e_{i_{1},a_{i_{1}}-1},0)=-1/4a_{i_{1}}^{2}$.
\qed
\end{pf}

The following Corollary~\ref{lem4-2} requires the condition {\rm (iv)} in Theorem~\ref{first}
for the case that $A=(2,2,2)$.
\begin{cor}\label{lem4-2}
If $a_{i}=2$,
then one has  
\[
c(4e_{i,1},0)=-\frac{1}{96}. 
\]
\end{cor}

\begin{pf}
We have $c(\gamma,0)=0$ if $\gamma-e_{i_1,j_1}-e_{i_2,j_2}\ge 0$ for $i_1\ne i_2$
by Proposition \ref{sep} (by the condition {\rm (iv)} for $A=(2,2,2)$).
Then, taking the coefficient in front of $t_{1,1}t_{2,1}t_{3,1}e^{t_{\mu_A}}$ in $WDVV((i,1),(i,1),\mu_A ,\mu_A)$,
one has
\[
2c(e_{1}+2e_{i,1},0)\cdot c(e_{1,1}+e_{2,1}+e_{3,1},1)+24c(2e_{i,1}+2e_{i,a_{i}-1},0)\cdot 2\cdot c(e_{1,1}+e_{2,1}+e_{3,1},1)=0.
\]
By the conditions {\rm (ii)} and {\rm (vi)} in Theorem~\ref{first},
$c(e_{1}+2e_{i,1},0)=1/4$ and $c(e_{1,1}+e_{2,1}+e_{3,1},1)=1$.
Hence one has $c(4e_{i,1},0)=-1/96$.
\qed
\end{pf}

\begin{prop}\label{0reconst}
Assume that $c(\alpha,0)$ and $c(\alpha,1)$ are reconstructed for $|\alpha|\le k+3$, $k\in\ZZ_{\ge 0}$.
Then $c(\gamma,0)$ and $c(\gamma,1)$ with $|\gamma |\le k+4$ are reconstructed
from $c(\alpha,0)$ and $c(\alpha,1)$ with $|\alpha|\le k+3$.
\end{prop}

\begin{pf}
We shall split the proof of Proposition \ref{0reconst} into following four steps. 

\begin{lem}[Step 1]\label{lem1-0reconst}
Let $\beta \in \ZZ^{\mu_{A}-2}$ be a non-negative element such that $|\beta|=k+1$.
Then $c(\beta +e_{i,j}+e_{i,j'}+e_{i,a_{i}-1},0)$ can be reconstructed from 
$c(\alpha,0)$ and $c(\alpha,1)$ with $|\alpha|\le k+3$.
\end{lem}

\begin{pf}
Without loss of generality, we can assume $i=1$.
First, we show that $c(\beta +e_{1,1}+e_{1,j+j'-1}+e_{1,a_{1}-1},0)$ can be reconstructed 
from $c(\alpha,0)$ and $c(\alpha,1)$ with $|\alpha|\le k+3$. 
One has $\deg (t^{\beta}t_{1,j+j'-1})=1$. By Proposition \ref{sep}, there exist $e_{1,l},e_{1,l'}$ such that
\begin{itemize}
\item $\beta +e_{1,j+j'-1}-e_{1,l}-e_{1,l'}\geq 0$, 
\item $\deg(t_{1,l})+\deg(t_{1,l'})\leq 1$.
\end{itemize}
We put $\beta' :=\beta +e_{1,1}+e_{1,j+j'-1}-e_{1,l}-e_{1,l'}$.
Taking the coefficient in front of $t^{\beta'}t_{2,1}t_{3,1}e^{t_{\mu_A}}$ in $WDVV((1,l),(1,l'),\mu_A,\mu_A)$,
one has
\begin{eqnarray*}
(\beta'_{1,l}+1)(\beta'_{1,l'}+1)(\beta'_{1,a_{1}-1}+1)\cdot c(\beta +e_{1,1}+e_{1,j+j'-1}+e_{1,a_{1}-1},0)
\cdot a_{1}\cdot c(e_{1,1}+e_{2,1}+e_{3,1},1)\\
+(known \ \ terms)=0.
\end{eqnarray*}
By the condition {\rm (vi)} in Theorem~\ref{first}, 
$c(\beta +e_{1,1}+e_{1,j+j'-1}+e_{1,a_{1}-1},0)$ can be reconstructed from  
$c(\alpha,0)$ and $c(\alpha,1)$ with $|\alpha|\le k+3$.

Next, we show that $c(\beta +e_{2,1}+e_{3,1}+e_{1,j+j'},1)$ can be reconstructed  
from $c(\alpha,0)$ and $c(\alpha,1)$ with $|\alpha|\le k+3$.
Taking the coefficient in front of  $t^{\beta }e^{t_{\mu_A }}$ in $WDVV((1,1),(1,j+j'-1),(2.1),(3.1))$, 
one has 
\begin{eqnarray*}
\lefteqn{s_{1,j+j'-1,a_{1}-j-j'}\cdot c(e_{1,1}+e_{1,j+j'-1}+e_{1,a_{1}-j-j'},0)\cdot a_{1}\cdot }\\
&&(\beta_{1,j+j'}+1)(\beta_{2,1}+1)(\beta_{3,1}+1)\cdot 
c(\beta +e_{2,1}+e_{3,1}+e_{1,j+j'},1)\\
&&+ (known \ \ terms)\\
&&+(\beta_{1,1}+1)(\beta_{1,j+j'-1}+1)(\beta_{1,a_{1}-1}+1)\cdot c(\beta +e_{1,1}+e_{1,j+j'-1}+e_{1,a_{1}-1},0)
\cdot a_{1}\cdot \\
&&c(e_{1,1}+e_{2,1}+e_{3,1},1)=0.
\end{eqnarray*}
By the previous argument and Proposition \ref{lem3},  $c(\beta +e_{2,1}+e_{3,1}+e_{1,j+j'},1)$ can be reconstructed from
$c(\alpha,0)$ and $c(\alpha,1)$ with $|\alpha|\le k+3$.

Finally, we show that $c(\beta +e_{1,j}+e_{1,j'}+e_{1,a_{1}-1},0)$ can be reconstructed  
from $c(\alpha,0)$ and $c(\alpha,1)$ with $|\alpha|\le k+3$.
Taking the coefficient in front of  $t^{\beta }t_{2,1}t_{3,1}e^{t_{\mu_A }}$ in $WDVV((1,j),(1,j'),\mu_A,\mu_A)$, 
one has
\begin{eqnarray*}
\lefteqn{\hspace{-33mm}{\rm (i)} \ (\beta_{1,j}+1)(\beta_{1,j'}+1)(\beta_{1,a_{1}-1}+1)\cdot
c(\beta +e_{1,j}+e_{1,j'}+e_{1,a_{1}-1},0)\cdot a_{1}\cdot c(e_{1,1}+e_{2,1}+e_{3,1},1)}\\ 
&&+ (known \ \ terms)\\ 
&&+s_{j,j'',a_{1}-j-j'}\cdot c(e_{1,j}+e_{1,j'}+e_{1,a_{1}-j-j'},0)\cdot \\
&&a_{1}\cdot (\beta_{1,j+j'}+1)\cdot c(\beta +e_{2,1}+e_{3,1}+e_{1,j+j'},1)=0\\ 
&&\text{if} \ \ a_{1}-j+a_{1}-j'\geq a_{1}+1,\\ 
\lefteqn{\hspace{-33mm}{\rm (ii)} \ (\beta_{1,j}+1)(\beta_{1,j'}+1)(\beta_{1,a_{1}-1}+1)\cdot
c(\beta +e_{1,j}+e_{1,j'}+e_{1,a_{1}-1},0)\cdot a_{1}\cdot c(e_{1,1}+e_{2,1}+e_{3,1},1)}\\ 
&&+ (known \ \ terms)\\ 
&&+c(e_{1,j}+e_{1,j'}+e_{1},0)\cdot 1\cdot c(\beta +e_{2,1}+e_{3,1},1)=0\\
&&\text{if} \ \ a_{1}-j+a_{1}-j'=a_{1}, \\
\lefteqn{\hspace{-33mm}{\rm (iii)} \ (\beta_{1,j}+1)(\beta_{1,j'}+1)(\beta_{1,a_{1}-1}+1)\cdot
c(\beta +e_{1,j}+e_{1,j'}+e_{1,a_{1}-1},0)\cdot a_{1}\cdot c(e_{1,1}+e_{2,1}+e_{3,1},1)}\\ 
&&+ (known \ \ terms)=0\\ 
&&\text{if} \ \ a_{1}-j+a_{1}-j'<a_{1}.
\end{eqnarray*}
By the second argument and Proposition \ref{lem3}, 
$c(\beta +e_{1,j}+e_{1,j'}+e_{1,a_{1}-1},0)$ can be reconstructed from
$c(\alpha,0)$ and $c(\alpha,1)$ with $|\alpha|\le k+3$.
\qed
\end{pf}

\begin{lem}[Step 2]\label{lem2-0reconst}
For a non-negative $\gamma \in \ZZ^{\mu_{A}-2}$ with $|\gamma|=k+4$,
$c(\gamma,1)$ can be reconstructed from $c(\alpha,0)$ and $c(\alpha,1)$ 
with $|\alpha|\le k+3$ 
\end{lem}

\begin{pf}
We shall split the proof of Lemma \ref{lem2-0reconst} into following three cases.  

\begin{sublem}[Step 2-Case 1]\label{sl1-l2-0reconst}
Let $i_{1}$, $i_{2}$, $i_{3}$ be distinct. Let 
$\gamma \in \ZZ^{\mu_{A}-2}$ be a non-negative element
such that $|\gamma|=k+4$ and
$\gamma-e_{i_{1},j_{1}}\geq 0$ and $2\leq j_{1}$.
Then $c(\gamma,1)$ can be reconstructed from $c(\alpha,0)$ and $c(\alpha,1)$ with $|\alpha|\le k+3$
\end{sublem}

\begin{pf}
We put $\gamma':=\gamma-e_{i_{1}.j_{1}}-e_{i_{2},1}-e_{i_{3},1}+e_{i_{1},1}+e_{i_{1},j_{1}-1}+e_{i_{1},a_{i_{1}}-1}$.
Taking the coefficient in front of  
$t^{\gamma-e_{i_{1},j_{1}}}e^{t_{\mu_A}}$
in $WDVV((i_{1},1),(i_{1},j_{1}-1),\mu_A,\mu_A)$, one has
\begin{eqnarray*}
\lefteqn{s_{1,j_{1}-1,a_{i_{1}}-j_{1}}\cdot c(e_{i_{1},1}+e_{i_{1},j_{1}-1}+e_{i_{1},a_{i_{1}}-j_{1}},0)\cdot a_{i_{1}} \cdot 
\gamma_{i_{1},j_{1}} \cdot c(\gamma ,1)}\\ 
&&+ (known\ \ terms)\\
&&+ \gamma'_{i_{1},1}\gamma'_{i_{1},j_{1}-1}\gamma'_{i_{1},a_{i_{1}}-1}\cdot c(\gamma',0)\cdot a_{i_{1}}\cdot c(e_{1,1}+e_{2,1}+e_{3,1},1)=0
\end{eqnarray*}
By Lemma \ref{lem1-0reconst}, 
$c(\gamma',0)$ can be reconstructed from $c(\alpha,0)$ and $c(\alpha,1)$ with $|\alpha|\le k+3$.
Hence $c(\gamma,1)$ can be reconstructed from $c(\alpha,0)$ and $c(\alpha,1)$ with $|\alpha|\le k+3$.
\qed
\end{pf}

\begin{sublem}[Step 2-Case 2]\label{sl2-l2-0reconst}
Let $\gamma \in \ZZ^{\mu_{A}-2}$ be a non-negative element
such that $|\gamma|=k+4$ and
$\gamma=\gamma_{1,1}e_{1,1}+\gamma_{2,1}e_{2,1}+\gamma_{3,1}e_{3,1}$ for $\gamma_{1,1}\gamma_{2,1}\gamma_{3,1}\ne 0$.
Then $c(\gamma,1)$ can be reconstructed from $c(\alpha,0)$ and $c(\alpha,1)$ with $|\alpha|\le k+3$.
\end{sublem}

\begin{pf}
By counting the degree of the term $t^{\gamma}e^{t_{\mu_{A}}}$, one has
\[
\deg(t^{\gamma}e^{t_{\mu_{A}}})>\deg(t_{1,1}t_{2,1}t_{3,1}e^{t_{\mu_{A}}})=2.
\]
Then one has $c(\gamma,1)=0$.
\qed
\end{pf}

\begin{sublem}[Step2-Case 3]\label{sl3-l2-0reconst}
Let $\gamma \in \ZZ^{\mu_{A-2}}$ be a non-negative element such that $|\gamma|=k+4$ and $\gamma =\gamma_{1,1}e_{1,1}+\gamma_{2,1}e_{2,1}+\gamma_{3,1}e_{3,1}$ for $\gamma_{1,1}\gamma_{2,1}\gamma_{3,1}=0$. 
Then $c(\gamma,1)$ can be reconstructed from $c(\alpha,0)$ and $c(\alpha,1)$ with 
$|\alpha|\le k+3$. 
\end{sublem}

\begin{pf}
Assume that $\gamma_{i_1,1}=0$.
We put $\gamma':=\gamma-e_{i_{2},1}-e_{i_{3},1}+e_{i_{1},1}+e_{i_{1},a_{i_{1}}-1}+e_{i_{1},a_{i_{1}}-1}$.
Taking the coefficient in front of $t_{1,1}^{\gamma_{1,1}}t_{2,1}^{\gamma_{2,1}}t_{3,1}^{\gamma_{3,1}}e^{t_{\mu_A}}$
in $WDVV((i_1,1),(i_1,a_{i_1}-1),\mu_A,\mu_A)$, one has 
\begin{eqnarray*}
\lefteqn{c(e_{i_1,1}+e_{1,a_{i_1}-1}+e_{1},0)\cdot c(\gamma ,1)
+(known \ \ terms)}\\
&&+\gamma'_{i_1,1}\gamma'_{i_1,a_{i_1}-1}\gamma'_{i_1,a_{i_1}-1}\cdot c(\gamma',0)\cdot a_{i_1}\cdot c(e_{1,1}+e_{2,1}+e_{3,1},1)=0
\end{eqnarray*}
One has $\gamma'-e_{i_2,1}\ge 0$ or $\gamma'-e_{i_3,1}\ge 0$.
Then one has $c(\gamma',0)=0$ by Proposition~\ref{sep}.
Hence $c(\gamma ,1)$ can be reconstructed from  
$c(\alpha,0)$ and $c(\alpha,1)$ with $|\alpha|\le k+3$.
Therefore $c(\gamma,1)$ can be reconstructed from  
$c(\alpha,0)$ and $c(\alpha,1)$ with $|\alpha|\le k+3$.
\qed
\end{pf}
Then we have Lemma~\ref{lem2-0reconst}.
\qed
\end{pf}

\begin{lem}\label{lem4.1}
If $a_{i_{1}}\geq 3$ and $i_{1},i_{2},i_{3}$ are distinct, 
then one has  
\[
c(e_{i_{1},j_{1}+1}+e_{i_{1},a_{i_{1}}-j_{1}}+e_{i_{2},1}+e_{i_{3},1},1)=
\begin{cases}
\frac{1}{a_{i_{1}}} \ \ \text{if} \ \ a_{i_{1}}-j_{1}\neq j_{1}+1,\\
\frac{1}{2a_{i_{1}}} \ \ \text{if} \ \ a_{i_{1}}-j_{1}=j_{1}+1.
\end{cases}
\]
\end{lem}

\begin{pf}
Taking the term in front of $t_{i_{2},1}t_{i_{3},1}e^{t_{\mu_A}}$ 
in $WDVV((i_{1},a_{i_{1}}-j_{1}-1),(i_{1},1),(i_{1},j_{1}+1),\mu_A)$, 
one has
\begin{eqnarray*}
\lefteqn{\hspace{-30mm}{\rm (i)} \ \frac{1}{a_{i_{1}}}\cdot a_{i_{1}} \cdot c(e_{i_{1},j_{1}+1}+e_{i_{1},a_{i_{1}}-j_{1}}+e_{i_{2},1}+e_{i_{3},1},1)}\\
&&-1\cdot 1\cdot c(e_{1}+e_{i_{1},a_{i_{1}}-j_{1}-1}+e_{i_{1},j_{1}+1},0)=0\\
&&\text{if} \ \ a_{i_{1}}-j_{1}\neq j_{1}+1,\ a_{i_{1}}-j_{1}-1\neq j_{1}+1,\\
\lefteqn{\hspace{-30mm}{\rm (ii)} \ \frac{1}{a_{i_{1}}}\cdot a_{i_{1}} \cdot c(e_{i_{1},j_{1}+1}+e_{i_{1},a_{i_{1}}-j_{1}}+e_{i_{2},1}+e_{i_{3},1},1)}\\
&&-1\cdot 1\cdot 2c(e_{1}+e_{i,a_{i}-i_{1}-1}+e_{i,i_{1}+1},0)=0\\
&&\text{if} \ \ a_{i_{1}}-j_{1}\neq j_{1}+1, \ a_{i_{1}}-j_{1}-1=j_{1}+1,\\
\lefteqn{\hspace{-30mm}{\rm (iii)} \ \frac{1}{a_{i_{1}}}\cdot a_{i_{1}} \cdot 2c(e_{i_{1},j_{1}+1}+e_{i_{1},a_{i_{1}}-j_{1}}+e_{i_{2},1}+e_{i_{3},1},1)}\\
&&-1\cdot 1\cdot c(e_{1}+e_{i_{1},a_{i_{1}}-j_{1}-1}+e_{i_{1},j_{1}+1},0)=0\\
&&\text{if} \ \ a_{i_{1}}-j_{1}=j_{1}+1.
\end{eqnarray*}
Hence one has Lemma \ref{lem4.1}.
\qed
\end{pf}

\begin{lem}[Step 3]\label{lem3-0reconst}
Let $\gamma \in \ZZ^{\mu_{A}-2}$ be a non-negative element such that 
$|\gamma|=k+4$ and $\gamma-e_{i,1}\geq 0$ for some $i$.
Then $c(\gamma, 0)$ can be reconstructed from $c(\alpha,0)$ and $c(\alpha,1)$ with $|\alpha|\le k+3$.
\end{lem}

\begin{pf}
Let $i_{1},i_{2},i_{3}$ be distinct. 
We will show this claim by the induction on the degree of parameter $t_{i,j}$.
By Lemma \ref{lem1-0reconst}, $c(\beta +e_{i_1,j}+e_{i_1,j'}+e_{i_1,a_{i_1}-1},0)$ 
with $|\beta|=k+1$ can be reconstructed from 
$c(\alpha,0)$ and $c(\alpha,1)$ with $|\alpha|\le k+3$.
Assume that $c(\gamma',0)$ with $|\gamma'|=k+4$ is known if $\gamma'-e_{i_1,1}-e_{i_1,n}\geq 0$, $n\ge l$.

We shall show that $c(\gamma,0)$ can be reconstructed from  
$c(\alpha,0)$ and $c(\alpha,1)$ with $|\alpha|\le k+3$
if $|\gamma|=k+4$ and $\gamma-e_{i_1,1}-e_{i_1,l-1}\geq 0$.
One has $\deg(t^{\gamma-e_{i_1,1}-e_{i_1,l-1}})=l/a_{i_1}$. 
By Proposition \ref{sep}, there exist $e_{i_1,m},e_{i_1,m'}$ such that 
\begin{itemize}
\item $\gamma-e_{i_1,1}-e_{i_1,l-1}-e_{i_1,m}-e_{i_1,m'}\geq 0$,
\item $\deg(t_{i_1,m})+\deg(t_{i_1,m'})\le l/a_{i_1}$.
\end{itemize}
Note that $\deg(t_{i_1,m})+\deg(t_{i_1,l})<1$ and $c(\alpha,1)$ with $|\alpha|=k+4$
can be reconstructed by Lemma~\ref{lem2-0reconst}.
We put $\beta:=\gamma-e_{i_1,l-1}-e_{i_1,m}-e_{i_1,m'}$.
Taking the coefficient in front of $t^{\beta}t_{i_2,1}t_{i_3,1}e^{t_{\mu_A}}$ in $WDVV((i_1,m),(i_1,m'),(i_1,l),\mu_A)$, 
one has 
\[
\gamma_{i_1,m}\gamma_{i_1,m'}\gamma_{i_1,l-1}\cdot c(\gamma ,0)\cdot a_{i_1}\cdot c(e_{i_1,a_{i_1}+1-l}+e_{i_1,l}+e_{i_2,1}+e_{i_3,1},1)
+(known \ \ terms)=0,
\]
By Lemma \ref{lem4.1}, $c(\gamma ,0)$ can be reconstructed from  
$c(\alpha,0)$ and $c(\alpha,1)$ with $|\alpha|\le k+3$
if $\gamma-e_{i,1}\ge 0$. 
\qed
\end{pf}

\begin{lem}[Step 4]\label{lem4-0reconst}
$c(\gamma,0)$ with $|\gamma|=k+4$ can be reconstructed from $c(\alpha,0)$ and $c(\alpha,1)$ with $|\alpha|\le k+3$.
\end{lem}

\begin{pf}
We will show this claim by the induction on the degree of parameter $t_{i,j}$.
By Lemma \ref{lem3-0reconst}, $c(\gamma ,0)$ can be reconstructed from 
$c(\alpha,0)$ and $c(\alpha,1)$ with $|\alpha|\le k+3$
if $\gamma-e_{i,1}\ge 0$. 
Assume that $c(\gamma',0)$ with $|\gamma'|=k+4$ is known if $\gamma'-e_{i,n}\geq 0$ and $n\leq l$.
We shall show $c(\gamma,0)$ can be reconstructed from $c(\alpha,0)$ and $c(\alpha,1)$ with $|\alpha|\le k+3$
if $\gamma-e_{i,l+1}\ge 0$.

Taking the coefficient in front of $t^{\gamma-e_{i,j}-e_{i,j'}-e_{i,l+1}}$ in $WDVV((i,1),(i,l),(i,j),(i,j'))$, one has 
\[
s_{1,l,a_{i}-1-l}\cdot c(e_{i,1}+e_{i,l}+e_{i,a_{i}-1-l},0)\cdot a_{i}\cdot \gamma_{i,j}\gamma_{i,j'}\gamma_{i,l+1}\cdot c(\gamma ,0)
+(known \ \ terms)=0.
\]
Hence $c(\gamma ,0)$ can be reconstructed from $c(\alpha,0)$ and $c(\alpha,1)$ with $|\alpha|\le k+3$
\qed
\end{pf}
Therefore we have Proposition~\ref{0reconst}
\qed
\end{pf}

By Proposition \ref{lem3}, Proposition \ref{lem3.1} and Proposition \ref{0reconst}, 
$c(\gamma,0)$ and $c(\gamma,1)$
can be reconstructed from $c(\beta,0)$ with $|\beta|=3$. 

\subsection{$c(\alpha,m)$ can be reconstructed}\label{reconst-s2}

In Sub-section \ref{reconst-s1}, we showed that $c(\alpha,0)$ and $c(\alpha,1)$ can be reconstructed from $c(\beta,0)$ with $|\beta|=3$.
We define the well-order $\prec$ on $\ZZ_{\ge 0}^{2}$ as follows:
\begin{itemize}
\item $(|\alpha|,m)\prec$ $(|\beta|,n)$ if $m<n$. 
\item $(|\alpha|,m)\prec$ $(|\beta|,m)$ if $|\alpha| <|\beta|$.
\end{itemize} 
We shall prove that $c(\alpha, m)$ can be reconstructed 
from $c(\beta,0)$ with $|\beta|=3$ by the induction on the well order $\prec$
on $\ZZ_{\ge 0}^{2}$.

\begin{prop}\label{mreconst}
$c(\gamma,m)$ with $m\ge 2$ can be reconstructed from $c(\beta,0)$ with $|\beta|=3$.
\end{prop}

\begin{pf}
Assume that $c(\alpha, n)$ can be reconstructed from $c(\beta,0)$, $|\beta|=3$ 
if $(|\alpha|,n)\prec (0,m-1)$.
First, we shall show that $c(0,m)$ can be reconstructed.

\begin{lem}\label{lem0-mreconst}
$c(0,m)$ can be reconstructed from $c(\alpha, n)$
with $(|\alpha|,n)\prec (0,m)$.
\end{lem}

\begin{pf}
Taking the coefficient in front of  
$e^{mt_{\mu_A}}$
in $WDVV((i,1),(i,a_{i}-1),\mu_A,\mu_A)$, one has
\[
c(e_{i,1}+e_{i,a_{i}-1}+e_{1},0)\cdot 1 \cdot m^{3}\cdot c(0,m)
+ (known\ \ terms)=0.
\]
Hence $c(0,m)$ can be reconstructed from $c(\alpha, n)$
with $(|\alpha|,n)\prec (0,m)$.
\qed
\end{pf}

Next, we shall split the second step of the induction into following three cases. 

\begin{lem}[Case 1]\label{lem1-mreconst}
Let $\gamma \in \ZZ^{\mu_{A}-2}$ be a non-negative element such that $|\gamma|=k+1$
and $\gamma-e_{i,j}\ge 0$ for $2\leq j$.
Then $c(\gamma,m)$ can be reconstructed from $c(\alpha, n)$
with $(|\alpha|,n)\prec (k+1,m)$.
 \end{lem}

\begin{pf}
Taking the coefficient in front of  
$t^{\gamma-e_{i,j}}e^{mt_{\mu_A}}$
in $WDVV((i,1),(i,j-1),\mu_A,\mu_A)$, one has
\[
s_{1,j-1,a_{i}-j}\cdot c(e_{i,1}+e_{i,j-1}+e_{i,a_{i}-j},0)\cdot a_{i} \cdot m^{2}\cdot \gamma_{i,j} \cdot c(\gamma ,m)
+ (known\ \ terms)=0,
\]
Hence $c(\gamma ,m)$ can be reconstructed from $c(\alpha, n)$
with $(|\alpha|,n)\prec (k+1,m)$.
\qed
\end{pf}

\begin{lem}[Case 2]\label{lem2-mreconst}
Let $\gamma \in \ZZ^{\mu_{A}-2}$ be a non-negative element such that $|\gamma|=k+1$
and $\gamma :=\gamma_{1,1}e_{1,1}+\gamma_{2,1}e_{2,1}+\gamma_{3,1}e_{3,1}$ for $\gamma_{1,1}\gamma_{2,1}\gamma_{3,1}\ne 0$.
Then $c(\gamma ,m)$ can be reconstructed from $c(\alpha, n)$
with $(|\alpha|,n)\prec (k+1,m)$.
\end{lem}

\begin{pf}
Taking the coefficient in front of $t_{1,1}^{\gamma_{1,1}}t_{2,1}^{\gamma_{2,1}}t_{3,1}^{\gamma_{3,1}}e^{mt_{\mu_A}}$
in $WDVV((i,1),(i,a_{i}-1),\mu_A,\mu_A)$, one has
\begin{eqnarray*}
\lefteqn{\hspace{-50mm}{\rm (i)} \ \{ c(e_{i,1}+e_{i,a_{i}-1}+e_{1},0)\cdot m^{3} +4\cdot 
c(2e_{i,1}+2e_{i,a_{i}-1},0)\cdot a_{i} 
\cdot m^{2}\cdot \gamma_{i,1}\} \cdot c(\gamma ,1)}\\ 
&&+ (known \ \ terms)=0 \\
&&\text{if} \ \ a_{i}\ge 3,\\
\lefteqn{\hspace{-50mm}{\rm (ii)} \ \{ 2c(2e_{i,1}+e_{1},0)\cdot m^{3} +24\cdot c(4e_{i,1},0)\cdot 2\cdot m^{2} 
\cdot \gamma_{i,1}\} \cdot c(\gamma ,1)}\\ 
&&+ (known \ \ terms)=0 \\
&&\text{if} \ \ a_{i}=2.
\end{eqnarray*}
If $\gamma_{i,1}\neq m$ for some $i$, $c(\gamma,m)$ can be reconstructed from $c(\alpha, n)$
with $(\alpha,n)\prec (k+1,m)$.
If $\gamma_{i,1}=m$ for all $i$, $\deg (t_{1,1}^{\gamma_{1,1}}t_{2,1}^{\gamma_{2,1}}t_{3,1}^{\gamma_{3,1}}e^{mt_{\mu_A}})$ 
is greater than $2$ except for the case $m=1$. 
\qed
\end{pf}

\begin{lem}[Case 3]\label{lem3-mreconst}
Let $\gamma \in \ZZ^{\mu_{A}-2}$ be a non-negative element such that $|\gamma|=k+1$ and $\gamma :=\gamma_{1,1}e_{1,1}+\gamma_{2,1}e_{2,1}+\gamma_{3,1}e_{3,1}$
for $\gamma_{1,1}\gamma_{2,1}\gamma_{3,1}=0$.
Then $c(\gamma ,m)$ can be reconstructed from $c(\alpha, n)$
with $(|\alpha|,n)\prec (k+1,m)$.
\end{lem}

\begin{pf}
Assume that $\gamma_{i,1}=0$.
Taking the coefficient in front of $t_{1,1}^{\gamma_{1,1}}t_{2,1}^{\gamma_{2,1}}t_{3,1}^{\gamma_{3,1}}e^{mt_{\mu_A}}$
in $WDVV((i,1),(i,a_{i}-1),\mu_A,\mu_A)$, one has
\[
c(e_{1}+e_{i,a_{i}-1}+e_{i,1},0)\cdot m^{3} \cdot c(\gamma ,m)
+(known \ \ terms)=0.
\]
Hence $c(\gamma ,m)$ can be reconstructed from $c(\alpha, n)$
with $(|\alpha|,n)\prec (k+1,m)$.
\qed
\end{pf}
Therefore we have Proposition~\ref{mreconst}.
\qed
\end{pf}
We finish the proof of Theorem~\ref{first}.

\section{The Gromov-Witten Theory for Orbifold Projective Lines} 

Let $\PP^1_A$ be the orbifold projective line with three orbifold points whose orders are 
$a_1, a_2, a_3$ respectively.
Denote by $\CC[[H^*_{orb}(\PP^1_A,\CC)]]$ the completed symmetric algebra of the dual space of $H^*_{orb}(\PP^1_A,\CC)$ 
and consider a formal manifold $M$ whose structure sheaf $\O_M$ and tangent sheaf $\T_M$ are given by 
the algebra $\CC[[H^*_{orb}(\PP^1_A,\CC)]]$ and $\T_M:=H^*_{orb}(\PP^1_A,\CC)\otimes_\CC\CC[[H^*_{orb}(\PP^1_A,\CC)]]$. 
The Gromov--Witten theory for orbifolds developed by Abramovich--Graber--Vistoli \cite{agv:1} and Chen--Ruan \cite{cr:1} 
gives us the following.
\begin{prop}[\cite{agv:1, cr:1}]
There exists a structure of a formal Frobenius manifold of rank $\mu_A$ and dimension one 
on $M$ whose non-degenerate symmetric $\O_M$-bilinear form $\eta$ on $\T_M$ is given by the Poincar\'{e} pairing.
\end{prop}
\begin{pf}
See Theorem 6.2.1 of \cite{agv:1} and Theorem 3.4.3 of \cite{cr:1}.
\qed
\end{pf}

The following theorem is the main result in this section:
\begin{thm}\label{second}
The conditions in Theorem~\ref{first} are satisfied by the Frobenius structure constructed
from the Gromov--Witten theory for $\PP^1_A$. 
\end{thm}
We shall check the conditions in Theorem~\ref{first} one by one.
\subsection{Conditions (i) and (ii)}\label{con(1,2)}
The orbifold cohomology group of $\PP^{1}_A$ is, as a vector space, just the
singular cohomology group of the inertia orbifold: 
\[
\I\PP^{1}_A=\PP^{1}_A \bigsqcup B(\ZZ/a_1 \ZZ)\bigsqcup B(\ZZ/a_2 \ZZ)\bigsqcup B(\ZZ/a_3 \ZZ),
\]
and the orbifold Poincar\'{e} pairing is given by twisting the usual Poincar\'{e} pairing:
\[
\displaystyle
\int_{\PP^{1}_A} \alpha \cup_{orb} \beta := \int_{\I\PP^{1}_A} \alpha \cup I \beta,
\]
where $I$ is the involution defined in \cite{agv:1, cr:1}.
Then we have the following Lemma~\ref{gw1}:
\begin{lem}\label{gw1}
We can choose a $\QQ$-basis $1=\Delta_{1},\Delta_{1,1},\dots,\Delta_{i,j},\dots,\Delta_{3,a_3-1},\Delta_{\mu_A}$ of 
the orbifold cohomology group $H^{*}_{orb}(\PP^{1}_{A},\QQ)$
such that
\[
H^{0}_{orb}(\PP^{1}_{A},\QQ )\simeq \QQ \Delta_{1}, \ \  
\Delta_{i,j}\in H^{2\frac{j}{a_{i}}}_{orb}(\PP^{1}_{A},\QQ ), \ \
H^{2}_{orb}(\PP^{1}_{A},\QQ )\simeq \QQ \Delta_{\mu_{A}}
\]
and
\[
\displaystyle\int_{\PP^{1}_{A}} \Delta_{1} \cup_{orb} \Delta_{\mu_A} =1, \ \
\displaystyle\int_{\PP^{1}_{A}} \Delta_{i,j} \cup_{orb} \Delta_{k,l} =
\begin{cases}
&\frac{1}{a_{i}} \ \ \text{if} \ \ k=i, \ l=a_{i}-j \\
&0 \ \ \ \text{otherwise} .
\end{cases}
\]
\end{lem}
\begin{pf}
The decomposition of $H^{*}_{orb}(\PP^{1}_A)$ follows from the decomposition of the inertia orbifold
$\I\PP^{1}_A$. The latter assertion immediately follows from the definition of the orbifold Poincar\'{e} pairing.
\qed
\end{pf}
Denote by $t_{1},t_{1,1},\dots ,t_{i,j},\dots , t_{3,a_3-1},t_{\mu_A}$ the dual coordinates of the $\QQ$-basis 
$\Delta_{1}$, $\Delta_{1,1},\dots,\Delta_{i,j},\dots,\Delta_{3,a_3-1},\Delta_{\mu_A}$ of  $H^{*}_{orb}(\PP^{1}_{A},\QQ)$ in Lemma~\ref{gw1}. 
Then, it is easy to see that the condition {\rm (ii)} is satisfied.
It is also easy to show
that the unit vector field $e\in\T_M$ and the Euler vector field $E\in\T_M$ are given as
\[
e=\frac{\p}{\p t_1}, \ \ E=t_1\frac{\p}{\p t_1}+\sum_{i=1}^3\sum_{j=1}^{a_i-1}\frac{a_i-j}{a_i}t_{i,j}\frac{\p}{\p t_{i,j}}
+\chi_{A}\frac{\p}{\p t_{\mu_A}},
\]
which is the condition {\rm (i)}.
\subsection{Condition (iii)}
The condition {\rm (iii)} follows from the divisor axiom and 
the definition of the genus zero potential $\F^{\PP^{1}_A}_{0}$.

\subsection{Condition (iv)}\label{con(4)}
The condition {\rm (iv)} is satisfied since the image of degree zero orbifold map with marked points on orbifold points 
on the source must be one of orbifold points on the target $\PP^1_A$.

\subsection{Condition (v)}

The orbifold cup product is the specialization of the quantum product
at $t_{1}=t_{1,1}=\dots=t_{3,a_3-1}=e^{t_{\mu_A}}=0$. 
Therefore, it turns out that the orbifold cup product can be determined by the degree zero
three point Gromov-Witten invariants.
\begin{lem}\label{gw5}
There is a $\CC$-algebra isomorphism between the orbifold cohomology ring $H^{*}_{orb}(\PP^{1}_{A},\CC)$ and 
$\CC[x_1,x_2,x_3]/I_{A}$, $I_{A}:=(x_1x_2,x_2x_3,x_3x_1,a_1x_1^{a_1}-a_2x_2^{a_2},a_2x_2^{a_2}-a_3x_3^{a_3},a_3x_3^{a_3}-a_1x_1^{a_1})$, 
where $\Delta_{1,1},\Delta_{2,1},\Delta_{3,1}$ are mapped to $x_1,x_2,x_3$, respectively.
\end{lem}
\begin{pf}
Under the same notation in Sub-section \ref{con(1,2)},
the orbifold cup product is given as follows:
\[
\Delta_{\alpha}\cup_{orb}\Delta_{\beta}=\sum_{\delta}\left<\Delta_{\alpha},\Delta_{\beta},\Delta_{\gamma}\right>^{\PP^{1}_A}_{0,3,0}
\eta^{\gamma\delta} \Delta_{\delta},
\]
where we set $\eta^{\gamma\delta}$ as follows:
\[
\displaystyle(\eta^{\gamma\delta})=(\int_{\PP^{1}_{A}}\Delta_{\gamma}\cup_{orb}\Delta_{\delta})^{-1}.
\]
By the previous argument in Sub-section \ref{con(4)}, 
we have
\[
\Delta_{i_1,j_1}\cup_{orb} \Delta_{i_2,j_2}=0 \ \text{if} \ i_1\ne i_2.
\] 
By the formula
\begin{align*}
 \int_{\PP^{1}_A} \Delta_{i_1,j_1} \cup_{orb} \Delta_{i_1,j'_1} \cup_{orb} \Delta_{i_1,j''_1}
=&\  \frac{1}{|\ZZ/a_{i_1}\ZZ|} \int_{pt} ev^{*}_1 (\Delta_{i_1,j_1})\cup ev^{*}_2 (\Delta_{i_1,j'_1})\cup ev^{*}_3 (\Delta_{i_1,j''_1})\\
=&\ 
\begin{cases}
\frac{1}{a_{i_1}} \ \ \ \text{if} \ j_1 +j'_1 +j''_1=a_{i_1}, \\
0\ \ \ \text{otherwise},
\end{cases}
\end{align*}
we have
\[
\Delta_{i_{1},j_{1}}\cup_{orb} \Delta_{i_1,j'_1}=\Delta_{i_1,j_1+j'_1} \ \text{if} \ j_1+j'_1\le a_{i_1}-1,
\]
and hence
\[
\displaystyle\Delta_{i_1,1}^{p}:=\underbrace{\Delta_{i_1,1}\cup_{orb}\cdots \cup_{orb} \Delta_{i_1,1}}_{a_{i_1} \ times}=
\frac{1}{a_{i_1}}\Delta_{\mu_A}.
\]
Therefore we have Lemma~\ref{gw5}.
\qed
\end{pf}

\begin{lem}\label{gw6}
The term \[
\begin{cases}
e^{t_{\mu_A}}\quad \textit{if}\quad a_{1}=a_{2}=a_{3}=1,\\
t_{3,1}e^{t_{\mu_A}}\quad \textit{if}\quad 1=a_{1}=a_{2}<a_{3},\\
t_{2,1}t_{3,1}e^{t_{\mu_A}}\quad \textit{if}\quad 1=a_{1}<a_{2},\\
t_{1,1}t_{2,1}t_{3,1}e^{t_{\mu_A}}\quad \textit{if}\quad a_{1}\ge 2,
\end{cases}
\] occurs with the coefficient $1$ in the $\F^{\PP^{1}_A}_{0}$.
\end{lem}
\begin{pf}
This lemma follows from the fact that the Gromov-Witten invariant counts the number of
orbifold maps from $\PP^{1}_{A}$ to $\PP^{1}_{A}$ of degree $1$ fixing three marked $($orbifold$)$ points, 
which is exactly the identity map. 
\qed
\end{pf}


\end{document}